\documentclass{amsart}

\usepackage{amssymb}
\usepackage{eucal}
\usepackage[all]{xy}
\usepackage{longtable}

\newcommand{\ZZ}{\mathbb{Z}}
\newcommand{\FF}{\mathbb{F}}
\newcommand{\CC}{\mathbb{C}}
\newcommand{\QQ}{\mathbb{Q}}

\newcommand{\EE}{\mathbb{E}}
\newcommand{\GG}{\mathbb{G}}
\newcommand{\KK}{\mathbb{K}}
\newcommand{\LL}{\mathbb{L}}

\newcommand{\PP}{\mathbb{P}}
\newcommand{\TT}{\mathbb{T}}

\newcommand{\Gm}{\GG_{\mathrm{m}}}
\newcommand{\Ga}{\GG_{\mathrm{a}}}

\newcommand{\cA}{\mathcal{A}}
\newcommand{\cH}{\mathcal{H}}
\newcommand{\cP}{\mathcal{P}}
\newcommand{\cR}{\mathcal{R}}
\newcommand{\cT}{\mathcal{T}}

\newcommand{\bF}{\mathbf{F}}
\newcommand{\bff}{\mathbf{f}}
\newcommand{\bg}{\mathbf{g}}
\newcommand{\bh}{\mathbf{h}}
\newcommand{\bL}{{L}}
\newcommand{\bm}{\mathbf{m}}

\newcommand{\bn}{\mathbf{n}}
\newcommand{\bnu}{\boldsymbol{\nu}}
\newcommand{\bp}{\mathbf{p}}
\newcommand{\bq}{\mathbf{q}}
\newcommand{\bsigma}{\boldsymbol{\sigma}}
\newcommand{\bu}{\mathbf{u}}
\newcommand{\bv}{\mathbf{v}}
\newcommand{\bw}{\mathbf{w}}
\newcommand{\bx}{\mathbf{x}}

\newcommand{\tLambda}{\widetilde{\Lambda}}
\newcommand{\tmu}{\tilde{\mu}}

\newcommand{\tPsi}{\widetilde{\Psi}}
\newcommand{\tpi}{\widetilde{\pi}}
\newcommand{\tSigma}{\widetilde{\Sigma}}

\newcommand{\fa}{\mathfrak{a}}
\newcommand{\fb}{\mathfrak{b}}

\newcommand{\fC}{\mathfrak{C}}
\newcommand{\fm}{\mathfrak{m}}

\newcommand{\fp}{\mathfrak{p}}
\newcommand{\fq}{\mathfrak{q}}

\newcommand{\rB}{\mathrm{B}}

\newcommand{\dP}{P^{\dagger}}
\newcommand{\dQ}{Q^{\dagger}}

\newcommand{\oalpha}{\overline{\alpha}}
\newcommand{\oF}{\overline{F}}

\newcommand{\ok}{\overline{k}}
\newcommand{\oK}{\overline{K}}
\newcommand{\oL}{\overline{L}}

\newcommand{\oPhi}{\overline{\Phi}}
\newcommand{\opsi}{\overline{\psi}}

\newcommand{\sC}{\mathsf{C}}
\newcommand{\sM}{\mathsf{M}}
\newcommand{\sN}{\mathsf{N}}
\newcommand{\sP}{\mathsf{P}}

\newcommand{\sm}{\mathsf{m}}
\newcommand{\sn}{\mathsf{n}}

\newcommand{\uphi}{\underline{\phi}}

\newcommand{\one}{\mathbf{1}}

\DeclareMathOperator{\Aut}{Aut}

\DeclareMathOperator{\End}{End}
\DeclareMathOperator{\den}{den}

\DeclareMathOperator{\DR}{DR}

\DeclareMathOperator{\GL}{GL}
\DeclareMathOperator{\Hom}{Hom}

\DeclareMathOperator{\im}{im}
\DeclareMathOperator{\ord}{ord}

\DeclareMathOperator{\Mat}{Mat}
\DeclareMathOperator{\MT}{MT}

\DeclareMathOperator{\rank}{rk}

\DeclareMathOperator{\Sol}{Sol}
\DeclareMathOperator{\Spec}{Spec}

\DeclareMathOperator{\Supp}{Supp}
\newcommand{\Tens}{\textstyle{\bigotimes}}
\newcommand{\tr}{\textnormal{tr}}

\newcommand{\iso}{\stackrel{\sim}{\to}}

\newcommand{\absI}[1]{{\left\lvert #1 \right\rvert}_\infty}
\newcommand{\absT}[1]{\left\| #1 \right\|}
\newcommand{\laurent}[2]{#1(\!( #2 )\!)}
\newcommand{\power}[2]{#1[\![ #2 ]\!]}
\newcommand{\genpow}[2]{#1\langle\!\langle #2 \rangle\!\rangle}
\newcommand{\tate}[2]{#1\{ #2\}}

\newcommand{\trdeg}[1]{{\textnormal{tr.\ deg}_{#1}\ }}

\newcommand{\Module}[1]{\mathbf{Mod}(#1)}
\newcommand{\Rep}[2]{\mathbf{Rep}(#1,#2)}

\newcommand{\Vector}[1]{\mathbf{Vec}(#1)}

\newcommand{\kts}{\ok(t)[\bsigma,\bsigma^{-1}]}
\newcommand{\ktsp}{\ok[t;\bsigma]}

\newcommand{\oFqt}{\overline{\FF_q(t)}}
\newcommand{\okt}{\overline{k(t)}}
\newcommand{\oLL}{\overline{\LL}}
\newcommand{\Qbar}{\overline{\QQ}}

\newtheorem{theorem}[subsubsection]{Theorem}
\newtheorem{lemma}[subsubsection]{Lemma}
\newtheorem{corollary}[subsubsection]{Corollary}
\newtheorem{proposition}[subsubsection]{Proposition}

\theoremstyle{remark}

\numberwithin{equation}{subsubsection}

\begin{document}

%%%%%%%%%%%%%%%%%%%%%%%%%%%%%%%%%%%%%%%%%
%               Title, Etc.             %
%%%%%%%%%%%%%%%%%%%%%%%%%%%%%%%%%%%%%%%%%

\title[Tannakian duality for $t$-motives and Carlitz logarithms]{Tannakian duality
for Anderson-Drinfeld motives and algebraic independence of
Carlitz~logarithms}
\author{Matthew A. Papanikolas}
\address{Department of Mathematics\\
  Texas A{\&}M University\\
  College Station, TX 77843}
\email{map@math.tamu.edu}
%%%
\thanks{Research supported by NSF grant DMS-0340812 and NSA grant
MDA904-03-1-0019}
%%%

\date{June 29, 2007}

%%%%%%%%%%%%%%%%%%%%%%%%%%%%%%%
%          Abstract           %
%%%%%%%%%%%%%%%%%%%%%%%%%%%%%%%
\begin{abstract}
We develop a theory of Tannakian Galois groups for $t$-motives and
relate this to the theory of Frobenius semilinear difference
equations.  We show that the transcendence degree of the period
matrix associated to a given $t$-motive is equal to the dimension of
its Galois group.  Using this result we prove that Carlitz
logarithms of algebraic functions that are linearly independent over
the rational function field are algebraically independent.
\end{abstract}

%%%%%%%%%%%%%%%%%%%%%%%%%%%%%%%%%%%%%%%%%
%               Document Text           %
%%%%%%%%%%%%%%%%%%%%%%%%%%%%%%%%%%%%%%%%%

\maketitle

\tableofcontents

\section{Introduction}

\subsection{Periods of $t$-motives}

\subsubsection{Notation}
Let $\FF_q$ be the field of $q$ elements, where $q$ is a power of a
prime $p$.  Let $k := \FF_q(\theta)$, where $\theta$ is
transcendental over $\FF_q$, and define an absolute value
$\absI{\,\cdot\,}$ at the infinite place of $k$ so that
$\absI{\theta} = q$.  Let $k_\infty := \laurent{\FF_q}{1/\theta}$ be
the $\infty$-adic completion of $k$, let $\overline{k_\infty}$ be an
algebraic closure, let $\KK$ be the $\infty$-adic completion of
$\overline{k_\infty}$, and let $\ok$ be the algebraic closure of $k$
in $\KK$.

\subsubsection{Anderson $t$-motives}
Let $t$ be a variable over $\FF_q$ that is independent from
$\theta$, and let $\ktsp$ be the ring of polynomials in $t$ and
$\bsigma$ over $\ok$ subject to the relations
\[
  ct = tc,\quad  \bsigma t = t \bsigma,\quad \bsigma c = c^{1/q} \bsigma, \quad c
  \in \ok.
\]
An Anderson $t$-motive is a left $\ktsp$-module $\sM$ that is free
and finitely generated as both a left $\ok[t]$-module and as a left
$\ok[\bsigma]$-module and that satisfies $(t-\theta)^n\sM \subseteq
\bsigma \sM$ for all $n$ sufficiently large (see
\S\ref{SS:AndersontMotives}).  Anderson $t$-motives were originally
defined in \cite{abp04}, where they were called ``dual
$t$-motives.''

\subsubsection{Rigid analytic triviality}
We let $\TT := \tate{\KK}{t}$ be the Tate algebra of power series in
$\power{\KK}{t}$ that are convergent on the closed unit disk in
$\KK$, and let $\LL \subseteq \laurent{\KK}{t}$ be its fraction
field. Let $\EE$ be the subring of $\TT$ consisting of power series
that are everywhere convergent and whose coefficients lie in a
finite extension of $k_\infty$.  Finally, for a Laurent series $f =
\sum_i a_i t^i \in \laurent{\KK}{t}$ and an integer $n \in \ZZ$, we
set $\sigma^{-n}(f) := f^{(n)} := \sum_i a_i^{q^n} t^i$.

If $\sM$ is an Anderson $t$-motive and $\sm \in \Mat_{r \times
1}(\sM)$ has entries comprising a $\ok[t]$-basis of $\sM$, then
there is a matrix $\Phi \in \Mat_r(\ok[t])$ representing
multiplication by $\bsigma$ on $\sM$ so that
\[
  \bsigma \sm = \Phi \sm
\]
and $\det\Phi = c(t-\theta)^s$ for some $c \in \ok^{\times}$ and $s
\geq 1$. The Anderson $t$-motive is rigid analytically trivial (see
Proposition~\ref{P:ARATrivialization}) if there is a matrix $\Psi
\in \GL_r(\TT)$ so that
\[
  \Psi^{(-1)} = \Phi\Psi.
\]
It can be shown that the entries of $\Psi$ are in fact in $\EE$ (see
Proposition~\ref{P:TateIsEntire}).

\subsubsection{Connection with $t$-modules}
\label{SSS:tModuleConnection} The category of rigid analytically
trivial Anderson $t$-motives is equivalent to the category of
uniformizable abelian $t$-modules defined over $\ok$, as in
\cite{and86}.  For a given Anderson $t$-motive $\sM$ and associated
$t$-module $E$, there is an explicit connection
\[
  \textnormal{periods of $E$} \quad \longleftrightarrow \quad
  \textnormal{$\ok$-linear combinations of entries of
  $\Psi(\theta)^{-1}$}.
\]
The details of this relationship will be the subject of a future
paper with Anderson, but examples are already seen in \S\ref{SS:RAT}
for the Carlitz motive (see also S.~K.~Sinha \cite[\S 5.2]{sinha97}
for examples involving special values of the function field
$\Gamma$-function).

\subsubsection{Remarks on $t$-motive terminology}
G.~Anderson introduced $t$-motives in \cite{and86}.  Later in
\cite{abp04} dual $t$-motives, which had several technical
advantages, were introduced.  The algebraic properties of these two
types of $t$-motives are essentially the same, and the two
categories are anti-equivalent to each other.  In this paper we will
follow the dual $t$-motive point of view only, and throughout we
refer to them as Anderson $t$-motives.  In the following paragraph
we discuss a third type of $t$-motive, defined properly in
\S\ref{SS:AndersontMotives}, which are our primary objects of study.

\subsubsection{Tannakian category of $t$-motives}
\label{SSS:TannakiantMotivesIntro} In \S\ref{SS:AndersontMotives} we
show that the category of rigid analytically trivial Anderson
$t$-motives up to isogeny embeds as a full subcategory of a neutral
Tannakian category $\cT$ over $\FF_q(t)$.  Objects in $\cT$ are
called simply $t$-motives, and throughout the paper the term
``$t$-motive'' will refer exclusively to an object in $\cT$. In
particular, from this standpoint all $t$-motives are rigid
analytically trivial.  Also objects in $\cT$ do not necessarily come
from pure Anderson $t$-motives in the sense of~\cite{and86}, and so
$\cT$ is a mixed category.

By Tannakian duality, for each object $M$ in $\cT$, the Tannakian
subcategory $\cT_M$ generated by $M$ satisfies an equivalence of
categories
\[
  \cT_M \approx \Rep{\Gamma_M}{\FF_q(t)},
\]
where $\Rep{\Gamma_M}{\FF_q(t)}$ is the category of finite
dimensional representations over $\FF_q(t)$ of some algebraic
subgroup $\Gamma_M \subseteq \GL_r$ defined over $\FF_q(t)$ (see
\S\ref{SS:GaloisGroupstMotives}). The group $\Gamma_M$ is called the
Galois group of $M$.

It should be noted that R.~Pink \cite{pinkH} has defined a category
$\cH$ of mixed Hodge structures for function fields that is a
neutral Tannakian category over $\FF_q(t)$.  He showed that the
category of rigid analytically trivial Anderson $t$-motives that are
also ``mixed'' embeds as a full subcategory of $\cH$.  It would be
interesting to investigate the relationships among Pink's Hodge
structures, the $t$-motives defined in this paper, and their
associated Galois groups.  In the end our category of $t$-motives is
best suited for our transcendence applications, so we do not pursue
further here the connections with Pink's work.  See also D.~Goss
\cite{goss94} for additional comparisons between $t$-motives and
motives over $\QQ$.

The following is the main theorem of this paper (restated later as
Theorem~\ref{T:TrDegDimGalGrp}).

\begin{theorem} \label{T:TrDegDimGalGrpIntro}
Let $M$ be a $t$-motive, and let $\Gamma_M$ be its Galois group.
Suppose that $\Phi \in \GL_r(\ok(t)) \cap \Mat_r(\ok[t])$ represents
multiplication by $\bsigma$ on $M$ and that $\det \Phi =
c(t-\theta)^s$, $c \in \ok^\times$. Let $\Psi$ be a rigid analytic
trivialization of $\Phi$ in $\GL_r(\TT) \cap \Mat_r(\EE)$.  Finally,
let $L$ be the subfield of $\overline{k_\infty}$ generated over
$\ok$ by the entries of $\Psi(\theta)$.  Then
\[
  \trdeg{\ok} L = \dim \Gamma_M.
\]
\end{theorem}

\subsubsection{Grothendieck's conjecture}
In light of \S\ref{SSS:tModuleConnection}, the statement of
Theorem~\ref{T:TrDegDimGalGrpIntro} can be thought of as a function
field version of Grothendieck's conjecture on periods of algebraic
varieties.  For an abelian variety $A$ over $\Qbar$ of dimension
$d$, let $P$ be the period matrix of $A$ that represents an
isomorphism between $H^1(A(\CC),\QQ) \otimes_{\QQ} \CC$ and
$H^1_{\DR}(A/\CC)$, with basis defined over $\Qbar$. Grothendieck's
conjecture is that
\[
  \trdeg{\Qbar} \Qbar(P) = \dim \MT(A),
\]
where $\MT(A)$ is the Mumford-Tate group of $A$ and is an algebraic
subgroup of $\GL_{2d} \times \Gm$ over $\QQ$.  P.~Deligne
\cite[Cor.~I.6.4]{delmil82} has proved that the dimension of
$\MT(A)$ is an upper bound for the transcendence degree.
Conjecturally the Mumford-Tate group is isomorphic to the motivic
Galois group of the motive $h_1(A) \oplus \QQ(1)$ over $\QQ$.  More
generally Grothendieck's period conjecture states that if $X$ is a
smooth variety over $\Qbar$, then
\[
  \trdeg{\Qbar} \Qbar(P(X)) = \dim
  \Gamma_{X}^{\textrm{mot}},
\]
where $P(X)$ is the period matrix of $X$ and
$\Gamma_{X}^{\textrm{mot}}$ is the motivic Galois group of $X$ over
$\QQ$.  It should be pointed out that by work of C.~Bertolin
\cite{bert02} many standard transcendence conjectures over $\Qbar$,
such as Schanuel's conjecture, follow from expanded versions of
Grothendieck's period conjecture.

\subsection{Algebraic independence of Carlitz logarithms}
One application of Theorem~\ref{T:TrDegDimGalGrpIntro} is a
characterization of algebraic relations over $\ok$ of Carlitz
logarithms of algebraic numbers.

\subsubsection{Carlitz exponential}
\label{SSS:expC} The Carlitz exponential is the power series
\[
  \exp_C(z) := z + \sum_{i=1}^\infty
\frac{z^{q^i}}{(\theta^{q^i}-\theta)(\theta^{q^i}-\theta^q) \cdots
(\theta^{q^i} - \theta^{q^{i-1}})}.
\]
As is well known (see \cite[Ch.~3]{goss:FF}, \cite[\S 2.5]{thakur}),
the function defined by $\exp_C$ converges everywhere on $\KK$, is
$\FF_q$-linear, and has kernel $\FF_q[\theta]\,\tpi$, where
\[
  \tpi := \theta \sqrt[q-1]{-\theta} \prod_{i=1}^\infty \Bigl( 1-
  \theta^{1-q^i} \Bigr)^{-1} \in
  k_\infty(\sqrt[q-1]{-\theta})^\times.
\]
The Carlitz exponential also satisfies the functional equation
\[
  \exp_C(\theta z) = \theta\exp_C(z) + \exp_C(z)^q, \quad z \in \KK.
\]
Moreover, this functional equation induces an exact sequence of
$\FF_q[t]$-modules,
\[
  0 \to \FF_q[\theta]\,\tpi \to \KK \to \fC(\KK) \to 0,
\]
where $\fC(\KK)$ is the $\FF_q[t]$-module of $\KK$-valued points on
the Carlitz module $\fC$ (see \S\ref{SSS:CarlitzModule}) and where
$t$ acts by multiplication by $\theta$ on the first two terms. The
number $\tpi$ is called the Carlitz period.

\subsubsection{Carlitz logarithm} \label{SSS:logC}
The Carlitz logarithm is the inverse of $\exp_C(z)$,
\[
  \log_C(z) := z + \sum_{i=1}^\infty
\frac{z^{q^i}}{(\theta - \theta^q)(\theta - \theta^{q^2}) \cdots
(\theta - \theta^{q^i})},
\]
which as a function on $\KK$ converges for all $z \in \KK$ with
$\absI{z} < \absI{\theta}^{q/(q-1)}$.  The Carlitz logarithm is
$\FF_q$-linear and satisfies the functional equation
\[
  \theta\log_C(z) = \log_C(\theta z) + \log_C(z^q),
\]
for all $z \in \KK$ where all three terms converge.

\subsubsection{Linear forms in Carlitz logarithms}
We recall a theorem of J.~Yu.  Suppose $\lambda_1, \dots, \lambda_r
\in \KK$ satisfy $\exp_C(\lambda_i) \in \ok$ for each $i=1, \dots,
r$.  As in the previous section there are many potential $k$-linear
relations among $\lambda_1, \dots, \lambda_r$. However, Yu proved
that these are the only possible linear relations over $\ok$ in the
following function field analogue of Baker's theorem on linear forms
in logarithms.

\begin{theorem}[{Yu \cite[Thm.~4.3]{yu97}}]
Suppose $\lambda_1, \dots, \lambda_r \in \KK$ satisfy
$\exp_C(\lambda_i) \in \ok$ for $i=1, \dots, r$.  If
$\lambda_1, \dots, \lambda_r$ are linearly independent over $k$,
then the numbers $1, \lambda_1, \dots, \lambda_r$ are linearly
independent over $\ok$.
\end{theorem}

Yu's result is an application of his far reaching Theorem of the
Sub-$t$-module \cite[Thm.~0.1]{yu97}, which characterizes all
$\ok$-linear relations among logarithms of points in $\ok$ on
general $t$-modules. Transcendence results about the Carlitz periods
and Carlitz logarithms go back to Carlitz and Wade in the 1940's.
For detailed accounts of the history of transcendence results for
Drinfeld modules, including Yu's theorem, see W.~D.~Brownawell
\cite{brown98} and D.~S.~Thakur \cite[Ch.~10]{thakur}.

\subsubsection{Algebraic independence of Carlitz logarithms}
In characteristic~$0$, Baker's theorem on linear forms in natural
logarithms of algebraic numbers is best known.  In the situation of
Carlitz logarithms we use Theorem~\ref{T:TrDegDimGalGrpIntro} to
prove the following theorem (restated later as
Theorem~\ref{T:AlgIndCarlitzLogs}).

\begin{theorem} \label{T:AlgIndCarlitzLogsIntro}
Let $\lambda_1, \dots, \lambda_r \in \KK$ satisfy $\exp_C(\lambda_i)
\in \ok$ for each $i=1, \dots, r$.  If $\lambda_1, \dots, \lambda_r$
are linearly independent over $k$, then they are algebraically
independent over~$\ok$.
\end{theorem}

It should be noted that, using Mahler's method, L.~Denis
\cite{denis} has proved the special case of this theorem where
$\lambda_1, \dots, \lambda_r$ are restricted to values of $\log_C$
on elements of $\FF_q(\theta^{1/e})$, $e \geq 1$, of degree in
$\theta$ less than $q/(q-1)$.

\subsection{Methods of proof}

\subsubsection{$\sigma$-semilinear difference equations}
The category of $t$-motives is a certain full subcategory in the
category of left $\kts$-modules which are finite dimensional as
$\ok(t)$-vector spaces. To every $t$-motive $M$ one can associate a
matrix $\Phi \in \GL_r(\ok(t))$ representing multiplication by
$\sigma$ and a rigid analytic trivialization $\Psi \in \GL_r(\LL)$
so that $\Psi^{(-1)} = \Phi\Psi$.  Here recall that $\LL$ is the
fraction field of the Tate algebra $\TT$.  Thus the columns of
$\Psi$ satisfy a system of $\sigma$-semilinear difference equations
in the sense of \cite{vdpsing97}, where $\sigma = (f \mapsto
f^{(-1)}) : \LL \iso \LL$, and we develop the theory of such
equations in this context in \S\ref{S:Galois}. In spirit this theory
is close to the Galois theory of differential equations and
difference equations in characteristic~$0$ \cite{andre01},
\cite{beuk92}, \cite{del90}, \cite{magid}, \cite{vdp99},
\cite{vdpsing97}, \cite{vdpsing03}.

In \S\ref{S:Galois} we develop the Picard-Vessiot theory for certain
kinds of difference equations for $\sigma$ and construct their
difference Galois groups (see Theorem~\ref{T:GammaGroup}). However,
careful attention must be paid to the fact that the fixed field of
$\sigma$ in $\ok(t)$ is $\FF_q(t)$.  The Galois theory of difference
equations developed by M.~van der Put and M.~F.~Singer
\cite{vdpsing97} is quite useful here, but it does not completely
apply because they fundamentally use that the field of fixed
elements under the difference automorphism is algebraically closed.
On the one hand, because the fixed field of $\sigma$ in $\LL$ is
also $\FF_q(t)$, the Galois groups we construct are themselves
defined over $\FF_q(t)$. However, that $\FF_q(t)$ is not
algebraically closed nor even perfect presents several difficulties
because in general the $\FF_q(t)$-valued points of the Galois group
need not be dense and the group itself need not be a priori smooth.

\subsubsection{$t$-motives and difference Galois groups}
Given a $t$-motive $M$ of dimension $r$ over $\ok(t)$, the
difference Galois group $\Gamma$ is a subgroup of $\GL_r$ over
$\FF_q(t)$. Let $\Sigma$ be the $\ok(t)$-subalgebra of $\LL$
generated by the entries of $\Psi$ and $\det(\Psi)^{-1}$, and let
$\Lambda$ be its fraction field.  The field $\LL$ is naturally a
left $\kts$-module via the automorphism $\sigma$, and $\Sigma$ and
$\Lambda$ are both $\sigma$-invariant. Then
\[
  \Gamma(\FF_q(t)) \cong \Aut_\sigma(\Sigma/\ok(t)),
\]
where the right-hand side is the group of automorphisms of $\Sigma$
over $\ok(t)$ that commute with $\sigma$.  Moreover, this
identification is compatible with base extensions of $\FF_q(t)$ (see
\S\ref{SSS:sigmaautos}--\ref{SSS:sigmaBaseExtensions}).

We work out an explicit description of $\Gamma(\oFqt)$ in
\S\ref{SS:GaloisAction}, and, using crucially that $\LL$ is a
separable extension of $\ok(t)$ and that $\ok(t)$ is algebraically
closed in $\Lambda$, we show that $\Gamma$ has the following
properties:
\begin{itemize}
\item $\Gamma$ is smooth over $\FF_q(t)$
(Theorem~\ref{T:Smoothness}(b));
\item $\dim \Gamma = \trdeg{\ok(t)} \Lambda$,
(Theorem~\ref{T:Smoothness}(c));
\item The elements of $\Lambda$ fixed by $\Gamma(\oFqt)$ are
precisely $\ok(t)$ (Theorem~\ref{T:GammaFixedIsoK}).
\end{itemize}
These properties are essential for proving in
Theorem~\ref{T:GammaPsiIsGammaM} that
\[
  \Gamma \cong \Gamma_M,
\]
where $\Gamma_M$ is the Galois group associated to $M$ by Tannakian
duality.

\subsubsection{The proof of Theorem~\ref{T:TrDegDimGalGrpIntro}}
The primary vehicle for proving this theorem is a $\ok$-linear
independence criterion from \cite[Thm.~3.1.1]{abp04}.  It is stated
here in Theorem~\ref{T:ABPLinearIndependence}.  We apply this
criterion to the rigid analytic trivializations of tensor powers of
$M$ so as to compare the dimensions of the $\ok$-span of monomials
of the entries of $\Psi(\theta)$ of a given degree and the
$\ok(t)$-span of monomials in the entries of $\Psi$.  Ultimately we
show that
\[
  \trdeg{\ok} L = \trdeg{\ok(t)} \Lambda,
\]
the latter of which is the same as the dimension of $\Gamma_M$.

\subsubsection{Carlitz logarithms}
For $\alpha_1, \dots, \alpha_r \in \ok^{\times}$ with
$\absI{\alpha_i} < \absI{\theta}^{q/(q-1)}$ for $i=1, \dots,
r$, we define a $t$-motive $X$ so that the field generated over
$\ok$ by the entries of its rigid analytic trivialization $\Psi$
evaluated at $t=\theta$ is precisely
\[
  L = \ok(\Psi(\theta)) = \ok(\tpi,
  \log_C(\alpha_1),\dots, \log_C(\alpha_r)).
\]
Moreover, we show that arbitrary logarithms are $k$-linear
combinations of logarithms of this form in a precise way.  We
determine a set of defining equations of the Galois group $\Gamma_X$
of $X$ in Theorem~\ref{T:CarlitzLinearRelations} each of which is a
linear polynomial over $\ok(t)$.  These linear relations each
produce a $k$-linear relation on the logarithms and $\tpi$. We then
use Theorem~\ref{T:TrDegDimGalGrpIntro} to show that all algebraic
relations must arise from these relations.

\subsection{Acknowledgements}
The author thanks D.~Brownawell, L.~Denis, C.-Y.~Chang, D.~Goss,
L.-C.~Hsia, M.~van der Put, D.~Thakur, and J.~Yu for many helpful
discussions on the contents of this paper.  He further thanks the
National Center for Theoretical Sciences in Hsinchu, Taiwan, where
many of the results in this paper were proved.  The author
especially thanks G.~Anderson and N.~Ramachandran for their
indispensable advice throughout this project.  Finally the author
thanks the referee for several useful suggestions.

\section{Notation and preliminaries}

\subsection{Table of symbols}
\begin{center}
\begin{longtable}{p{0.66in}@{\hspace{5pt}$=$\hspace{5pt}}p{4in}}
$\power{R}{x}$ & power series in $x$ with
  coefficients in a ring $R$ \\
$\laurent{R}{x}$ & Laurent series in $x$ with
  coefficients in a ring $R$ \\
$\genpow{R}{x}$ & generalized power series ring in $x$ with
coefficients in a ring~$R$ \\
$\FF_q$ & finite field with $q = p^m$ elements \\
$k$ & $\FF_q(\theta) =$ rational functions in the variable $\theta$ over $\FF_q$ \\
$k_\infty$ & $\laurent{\FF_q}{1/\theta} = \infty$-adic
  completion of $k$ \\
$\overline{k_\infty}$ & algebraic closure of $k_\infty$ \\
$\zeta_{\theta}$ & a fixed $(q-1)$-th root of $-\theta$ in
$\overline{k_\infty}$ \\
$\KK$ & completion of $\overline{k_\infty}$ \\
$\ok$ & algebraic closure of $k$ in $\KK$ \\
$\TT$ & $\tate{\KK}{t} =$ ring of restricted power series; series
in $\power{\KK}{t}$ that converge on the closed unit disk
$\absI{t} \leq 1$ \\
$\LL$ & fraction field of $\TT$ \\
$M^{\vee}$ & dual vector space of a vector space $M$ \\
$\Vector{F}$ & category of finite dimensional vector spaces over a
field $F$ \\
$\Rep{\Gamma}{F}$ & for a field $F$ the category of finite
dimensional $F$-representations of an affine group scheme $\Gamma$
over $F$
\end{longtable}
\end{center}

\subsection{Preliminaries}

\subsubsection{Norms}
We let $\absI{\,\cdot\,}$ denote a fixed $\infty$-adic norm on
$\KK$. For a matrix $E \in \Mat_{r \times s}(\KK)$, we set $\absI{E}
= \sup \absI{E_{ij}}$.  For matrices $E$ and $F$, we observe that
$\absI{E+F} \leq \max( \absI{E},\absI{F})$ and $\absI{EF} \leq
\absI{E} \cdot \absI{F}$.

\subsubsection{Generalized power series}
Let $F$ be a field of characteristic $p$.  For a formal series $f :=
\sum_{i \in \QQ} a_i t^i$ with $a_i \in F$, we let $\Supp(f) := \{ i
\in \QQ \mid a_i \neq 0\}$.  We let $\genpow{F}{t}$ be the set of
such series for which $\Supp(f)$ is a well-ordered subset of $\QQ$.
This condition implies that $\genpow{F}{t}$ is a field under the
natural addition and multiplication of these series so that $t^it^j
= t^{i+j}$ (see P.~Ribenboim \cite[\S 2]{ribe92}).  If $F$ is
algebraically closed, then $\genpow{F}{t}$ is algebraically closed
\cite[\S 5]{ribe92}. If $F$ is a perfect field, then $\genpow{F}{t}$
is also perfect.

It should be noted that, when $F$ is algebraically closed,
$\genpow{F}{t}$ is not the algebraic closure of the Laurent series
field $\laurent{F}{t}$. For an explicit description of the field
$\overline{\laurent{F}{t}} \subseteq \genpow{F}{t}$, the reader is
directed to K.~Kedlaya~\cite{kedl01}.

By considering the inclusions
\[
  \FF_q(t) \subseteq \ok(t)
  \subseteq \laurent{\KK}{t} \subseteq \genpow{\KK}{t},
\]
we fix once and for all the inclusions of algebraically closed
fields
\[
  \overline{\FF_q(t)} \subseteq \overline{k(t)}
   \subseteq \overline{\laurent{\KK}{t}}
  \subseteq \genpow{\KK}{t}.
\]

\subsubsection{Entire functions}
A power series $f = \sum_{i=0}^\infty a_i t^i \in \power{\KK}{t}$ that
satisfies
\[
  \lim_{i \to \infty} {\textstyle \sqrt[i]{\absI{a_i}}} = 0
\]
and
\[
  [ k_\infty (a_0, a_1, a_2, \ldots) : k_\infty] < \infty,
\]
is an \emph{entire power series}.  As a function of $t$, such a
power series $f$ converges on all of $\KK$, and, when restricted to
$\overline{k_\infty}$, $f$ takes values in $\overline{k_\infty}$.
The ring of entire power series is denoted $\EE$.

\subsubsection{Restricted Laurent series}
A power series $\sum_{i=0}^\infty a_i t^i \in \power{\KK}{t}$ that
satisfies
\[
  \lim_{i \to \infty} \absI{a_i} = 0,
\]
is called a \emph{restricted power series}.  As functions of $t$,
these power series converge on the closed unit disk in $\KK$.  The
restricted power series form a subring $\TT=\tate{\KK}{t}$ of
$\power{\KK}{t}$, and $\EE$ is a subring of $\TT$.  The fraction
field of $\TT$, denoted $\LL$, is the field of \emph{restricted
Laurent series}.

Now at each point $a \in \KK$ with $\absI{a} \leq 1$, a function $f
\in \LL$ has a well-defined order of vanishing $\ord_a(f)$, and for
all but finitely many $\absI{a} \leq 1$, we have $\ord_a(f) = 0$.
Also each $f \in \LL$ has a unique factorization
\begin{equation} \label{E:LLfactorization}
  f = \lambda \biggl[\prod_{|a|_{\infty} \leq 1} (t - a)^{\ord_a(f)}
  \biggr] \biggl[1 + \sum_{i=1}^\infty b_i t^i \biggr],
\end{equation}
where $0 \neq \lambda \in \KK$, $\sup \absI{b_i} < 1$, and
$\absI{b_i} \to 0$ (see \cite[Cor.~2.2.4]{frevdp04}).  The series $1
+ \sum b_i t^i$ is a unit in $\TT$, and it follows that $\TT$ is a
principal ideal domain with maximal ideals generated by each $t-a$,
$\absI{a}\leq 1$ (see \cite[Thm.~2.2.9]{frevdp04}).

For $f = \sum_{i=0}^\infty a_i t^i \in \TT$, we define its norm
$\absT{f}$ to be
\[
  \absT{f} := \sup_i \absI{a_i} = \max_i \absI{a_i}.
\]
If $f \in \TT$ is written as in \eqref{E:LLfactorization}, then
$\absT{f} = \absI{\lambda}$. The norm $\absT{\,\cdot\,}$ is a
complete ultrametric norm on $\TT$ and satisfies
\begin{alignat*}{2}
  \absT{cf} &= \absI{c} \absT{f},\qquad &&\forall\, c \in \KK,
      f \in \TT, \\
  \absT{fg} &= \absT{f}\cdot \absT{g}, \qquad && \forall\, f, g \in \TT.
\end{alignat*}

\subsubsection{Twisting} \label{SSS:Twisting}
We define an automorphism $\sigma : \genpow{\KK}{t} \to
\genpow{\KK}{t}$ by setting
\[
 \sigma \biggl( \sum_{i \in \QQ} a_i t^i \biggr) :=
 \sum_{i \in \QQ} a_i^{1/q} t^i.
\]
If $f \in \genpow{\KK}{t}$ and $n \in \ZZ$, the \emph{$n$-fold
twist} of $f$ is defined to be
\[
  f^{(n)} := \sigma^{-n}(f).
\]
The automorphism $\sigma$ of $\genpow{\KK}{t}$ induces automorphisms
of several subrings, notably $\ok[t]$, $\ok(t)$, $\EE$, $\TT$,
$\LL$, $\power{\KK}{t}$, $\laurent{\KK}{t}$.  Moreover, $\sigma$
also leaves $\oFqt$, $\okt$, $\oLL$, and
$\overline{\laurent{\KK}{t}}$ invariant.

If $F$ is a subring of $\genpow{\KK}{t}$ that is invariant under
$\sigma$, we set
\[
  F^\sigma := \{ f \in F \mid \sigma(f) = f \}
\]
to be the elements of $F$ fixed by $\sigma$.  It is clear that
$F^{\sigma^n}$ is a subring of $F$ and that $F^{\sigma^m} \subseteq
F^{\sigma^n}$ if $m \mid n$.  For example,
\begin{alignat*}{2}
\genpow{\KK}{t}^\sigma &= \genpow{\FF_q}{t}, & \qquad \okt^\sigma &=
\overline{\FF_q(t)}^\sigma = \overline{\FF_q(t)} \cap
\genpow{\FF_q}{t}, \\
\laurent{\KK}{t}^\sigma &= \laurent{\FF_q}{t}, & \qquad
 \ok(t)^\sigma &= \FF_q(t).
\end{alignat*}
The only item that requires any explanation here is the description
of $\okt^\sigma$. For $\alpha \in \okt^\sigma$, let $x^m + b_{m-1}
x^{m-1} + \dots + b_0 \in \ok(t)[x]$ be the minimal polynomial of
$\alpha$ over $\ok(t)$. Since $\sigma(\alpha) = \alpha$, we have
that $\alpha$ is also a root of $x^m + \sigma(b_{m-1}) x^{m-1} +
\dots + \sigma(b_0)$. Taking the difference of these two relations,
we see that $\sigma(b_i)=b_i$ for each $i$, and so the minimal
polynomial of $\alpha$ has coefficients in $\ok(t)^\sigma =
\FF_q(t)$.

If $F$ is a matrix with entries in $\genpow{\KK}{t}$, then
$\sigma^{-n}(F) := F^{(n)}$ is defined by the rule
$\bigl(F^{(n)}\bigr)_{ij} := \bigl(F_{ij}^{(n)}\bigr)$.  If $F \in
\Mat_{r \times s}(\LL)$, we set $\absT{F} := \max_{i,j} \bigl\|
F_{ij} \bigr\|$, in which case $\bigl\| F^{(n)} \bigr\| =
\absT{F}^{q^n}$.

\begin{lemma} \label{L:ConstantTwist}
  For any $\alpha \in \KK$, there is a positive integer $s$ so that
  with respect to $\absI{\,\cdot\,}$ on $\KK$,
\[
  \lim_{n \to \infty} \alpha^{(ns)} = \begin{cases}
0 & \textnormal{if $\absI{\alpha} < 1$,} \\
c \in \overline{\FF}_q^{\times} & \textnormal{if $\absI{\alpha} = 1$,} \\
\infty & \textnormal{if $\absI{\alpha} > 1$.}
\end{cases}
\]
\end{lemma}

\begin{proof}
If $\absI{\alpha} \neq 1$, then the result is clear. Otherwise,
there is a unique $c \in \overline{\FF}_q^{\times}$ so that
$\absI{\alpha-c} < 1$. (See \cite[Lem.~2.4.4]{sinha97}.)  Then $c
\in \FF_{q^s}$ for some $s\geq 1$, and the result follows.
\end{proof}

\begin{lemma} \label{L:LLTwist}
  For any $f \in \TT$ with $\absT{f} \leq 1$, there is a positive
  integer $s$ so that with respect to $\absT{\,\cdot\,}$ on $\TT$,
\[
  \lim_{n \to \infty} f^{(ns)} \in \overline{\FF}_q[t].
\]
Also $\absT{f} = 1$ if and only if $\lim_{n \to \infty} f^{(ns)}
\neq 0$.
\end{lemma}

\begin{proof}
We use the factorization of $f$ in \eqref{E:LLfactorization}.  For
each $a$, with $\absI{a} \leq 1$, if $\ord_a(f) \neq 0$, then as in
Lemma~\ref{L:ConstantTwist} choose $s_a \geq 1$ and $c_a \in
\overline{\FF}_q$ so that $\lim_{n\to \infty} a^{(ns_a)} = c_a$.
Likewise, since $\absT{f} \leq 1$, we have $\absI{\lambda} \leq 1$,
and so we can choose $s_{\lambda} \geq 1$ and $c_{\lambda} \in
\overline{\FF}_q$ with $\lambda^{(ns_\lambda)} \to c_\lambda$.  Then
we let $s$ be the least common multiple of all the $s_a$'s and
$s_\lambda$.  From \eqref{E:LLfactorization}, with respect to
$\absT{\,\cdot\,}$,
\[
  \lim_{n \to \infty} \biggl[ 1 + \sum_{i=1}^\infty b_i t^i \biggr]^{(ns)}
  = 1,
\]
since $\sup \absI{b_i} < 1$.  Therefore,
\[
  \lim_{n \to \infty} f^{(ns)}
 = \lim_{n \to \infty} \lambda^{(ns)} \prod_{\absI{a} \leq 1}
   \bigl(t-a^{(ns)}\bigr)^{\ord_a(f)}
 = c_\lambda \prod_{\absI{a} \leq 1} (t-c_a)^{\ord_a(f)},
\]
which is in $\overline{\FF}_q[t]$.  Furthermore, $\absT{f} = 1$ if
and only if $\absI{\lambda} = 1$, which holds if and only if
$c_\lambda \neq 0$. Thus $\absT{f} = 1$ if and only if $\lim_{n \to
\infty} f^{(ns)} \neq 0$.
\end{proof}

\section{$t$-motives and Tannakian categories}

Here we will define a category $\cT$ of $t$-motives that is a
neutral Tannakian category over $\FF_q(t)$.  For all definitions of
tensor categories and Tannakian categories, we follow Deligne and
J.~S.~Milne \cite[\S II]{delmil82}.  Other useful references include
\cite{breen94}, \cite{del90}, \cite[App.~B]{vdpsing03}.

As mentioned in \S\ref{SSS:TannakiantMotivesIntro}, Tannakian
categories for $t$-motives have been considered previously by
Pink~\cite{pinkH}, though through a different construction. Parts of
the theory of $t$-motives defined below have been considered by
Y.~Taguchi \cite{tag95} and A.~Tamagawa \cite{tam95} in their study
of the Tate conjecture for $t$-modules. Also our theory has
similarities with the theory of $\sigma$-bundles defined by U.~Hartl
and Pink~\cite{harpin04}.

\subsection{The rings $\ktsp$ and $\kts$}

\subsubsection{Definition}
The ring $\kts$ is the noncommutative ring of Laurent polynomials in
the variable $\bsigma$ with coefficients in $\ok(t)$, subject to the
relation
\[
  \bsigma f = \sigma(f)\bsigma = f^{(-1)}\bsigma
\]
for all $f \in \ok(t)$.  Thus every element of $\kts$ has the form
$\sum_{i=-m}^m f_i \bsigma^i$, where $f_i \in \ok(t)$.

\subsubsection{Ring-theoretic properties}
The polynomials in $\bsigma$ with coefficients in $\ok[t]$ comprise
the subring $\ktsp$ of $\kts$.  The ring $\ok[\bsigma]$ is the
subring of polynomials with coefficients in $\ok$.   Both $\ktsp$
and $\kts$ are domains.  The center of $\ktsp$ is $\FF_q[t]$, and
the center of $\kts$ is $\FF_q(t)$.  The fundamental properties of
the ring $\ktsp$ are covered in \cite[\S 4]{abp04}.

\subsection{Pre-$t$-motives}

Here we define the category $\cP$ of pre-$t$-motives and explore its
basic properties.  In particular we show in
Theorem~\ref{T:PisRigidABTensor} that $\cP$ is a rigid abelian
$\FF_q(t)$-linear tensor category.

\subsubsection{The category $\cP$}
We let $\cP$ be the category of left $\kts$-modules that are finite
dimensional over $\ok(t)$.  Morphisms in $\cP$ are left
$\kts$-module homomorphisms.  We call $\cP$ the category of
\emph{pre-$t$-motives}, though it is worth noting that $\cP$ is a
category of difference modules with respect to the automorphism
$\sigma: \ok(t) \to \ok(t)$ in the sense of \cite{vdpsing97}.

\subsubsection{Preliminary properties of $\cP$}
The category of pre-$t$-motives is an abelian category.  For two
objects $P$ and $Q$ in $\cP$, it follows that $\Hom_{\cP}(P,Q)$ is
an $\FF_q(t)$-vector space.  A straightforward adaptation of the
proof of \cite[Thm.~2]{and86} shows that the map
\[
  \Hom_{\cP}(P,Q) \otimes_{\FF_q(t)} \ok(t) \to \Hom_{\ok(t)}(P,Q)
\]
is injective.  Thus $\Hom_{\cP}(P,Q)$ is a finite dimensional
$\FF_q(t)$-vector space.

\subsubsection{Representations of pre-$t$-motives}
\label{SSS:RepsPretmot} Given a $\ok(t)$-vector space $P$ and $p_1,
\dots, p_r \in P$, we call the vector
\[
\bp = \left[ \begin{matrix} p_1 \\ \vdots \\ p_r
\end{matrix} \right] \in \Mat_{r\times 1}(P)
\]
a \emph{basis for $P$} if $p_1, \dots, p_r$ form a $\ok(t)$-basis
for $P$.  If $P$ is a pre-$t$-motive, then there is a unique matrix
$\Phi = \Phi_{\bp} \in \GL_r(\ok(t))$ such that
\[
  \bsigma \bp = \Phi \bp.
\]
We say that \emph{$\Phi$ represents multiplication by $\bsigma$ on
$P$}.  Moreover, the matrix $\Phi \in \GL_r(\ok(t))$ uniquely
determines the left $\kts$-module structure on $P$ with respect to
$\bp$.

Now suppose that $\phi: P \to Q$ is a morphism of pre-$t$-motives
and that $\bp \in \Mat_{r\times 1}(P)$ and $\bq \in \Mat_{s\times
1}(Q)$ are bases for $P$ and $Q$ respectively.  If $B \in \Mat_{r
\times s}(\ok(t))$ represents $\phi$ as a map of $\ok(t)$-vector
spaces such that
\[
  \phi(\bff \cdot \bp) = \bff \cdot B\cdot\bq, \quad \bff \in \Mat_{1 \times
  r}(\ok(t)),
\]
then
\[
  B^{(-1)}\Phi_{\bq} = \Phi_{\bp} B.
\]
In particular, if $\bq$ is simply another basis of $P$, and $B \in
\GL_r(\ok(t))$ is the change of basis matrix, then $\Phi_{\bp} =
B^{(-1)}\Phi_{\bq}B^{-1}$.

\subsubsection{Tensor products of pre-$t$-motives}
Let $P$ and $Q$ be pre-$t$-motives.  Then the $\ok(t)$-vector
space $P \otimes_{\ok(t)} Q$ is made into a $\kts$-module by
defining
\[
  \bsigma(m \otimes n) := (\bsigma m) \otimes (\bsigma n).
\]
It is clear that then multiplication by $\bsigma$ is bijective on $P
\otimes_{\ok(t)} Q$ and that $P \otimes_{\ok(t)} Q$ is a
pre-$t$-motive.  Likewise we define arbitrary finite tensor products
of pre-$t$-motives with diagonal $\bsigma$-action.  For a fixed
pre-$t$-motive $P$ and $n \geq 1$, we set $P^{\otimes n} :=
\Tens_{i=1}^n P$ to be the $n$-th tensor power of $P$.

\subsubsection{Representations of tensor products}
Let $\bp = [p_1, \dots, p_r]^\tr$ and $\bq = [q_1, \dots,
q_s]^\tr$ be $\ok(t)$-bases for pre-$t$-motives $P$ and $Q$
respectively. Then, with respect to the basis
\[
\bp \otimes \bq := [ p_1 \otimes q_1, p_1 \otimes q_2, \dots, p_r
\otimes q_s ]^\tr,
\]
on $P \otimes Q$, the Kronecker product, $\Phi_{\bp \otimes \bq} =
\Phi_\bp \otimes \Phi_{\bq}$, represents multiplication by $\bsigma$
on $P \otimes Q$.  Similarly these conventions extend to arbitrary
finite tensor products of pre-$t$-motives.

\subsubsection{The Carlitz motive} \label{SSS:CarlitzMotive}
We define the \emph{Carlitz motive} to be the pre-$t$-motive $C$
whose underlying $\ok(t)$-vector space is $\ok(t)$ itself and on
which $\bsigma$ acts by
\[
  \bsigma f := (t-\theta)f^{(-1)}, \quad f \in C.
\]
For $n \geq 1$, the underlying $\ok(t)$-vector space of $C^{\otimes
n}$ is also $\ok(t)$, and multiplication by $\bsigma$ on $C^{\otimes
n}$ is given by
\[
  \bsigma f = (t-\theta)^n f^{(-1)}, \quad f \in C^{\otimes n}.
\]
See also \cite{andth90}.

\subsubsection{Internal Hom}
Let $P$ and $Q$ be pre-$t$-motives, and set
\[
  R := \Hom_{\ok(t)}(P,Q).
\]
Then $R$ is a $\ok(t)$-vector space.  We define a $\kts$-module
structure on $R$ by setting
\[
\bsigma\cdot \rho := \bsigma \circ \rho \circ \bsigma^{-1}, \quad
\rho \in R.
\]
It is straightforward to check that $\bsigma \cdot \rho : P \to Q$
is $\ok(t)$-linear, and so $\bsigma : R \to R$, and that this action
of $\bsigma$ extends naturally to a left $\kts$-module structure on
$R$. We write $\Hom(P,Q)$ for the $\kts$-module $R$ just defined. It
is also a pre-$t$-motive.

\subsubsection{Identity object} \label{SSS:one}
Let $\one := \ok(t)$ together with a $\bsigma$-action defined by
\[
  \bsigma f = \sigma(f) = f^{(-1)}, \quad f \in \one.
\]
Then $\one$ is a pre-$t$-motive.  Moreover, for any pre-$t$-motive
$P$, the natural isomorphisms, $P \otimes \one \cong \one \otimes P
\cong P$, are isomorphisms of pre-$t$-motives.  Thus $\one$ is an
identity object with respect to tensor products in $\cP$.

\begin{lemma} \label{L:End1}
$\End_{\cP}(\one) \cong \FF_q(t)$.
\end{lemma}

\begin{proof}
Suppose $\phi : \one \to \one$ is a morphism in $\cP$.  As a map of
$\ok(t)$-vector spaces, there is some $a \in \ok(t)$ so that
$\phi(f) = af$ for all $f \in \ok(t)$.  Since $\phi$ is also
$\kts$-linear, we must have $\bsigma a = a \bsigma$, and so $a$ is
in the center of $\kts$.  Thus $a \in \FF_q(t)$.
\end{proof}

\subsubsection{Duals}
Let $P$ be a pre-$t$-motive.  Then set
\[
  P^{\vee} := \Hom(P,\one).
\]
The pre-$t$-motive $P^{\vee}$ is called the \emph{dual of $P$}. As a
$\ok(t)$-vector space, $P^{\vee}$ is the dual vector space of $P$.
If $\bp$ forms a basis for $P$, let $\bp^{\vee}$ be the dual basis.
We find easily that
\[
  \Phi_{\bp^{\vee}} = \bigl(\Phi_{\bp}^{-1}\bigr)^{\tr}.
\]
If $\phi : P \to Q$ is a morphism of pre-$t$-motives, then the dual
morphism of $\ok(t)$-vector spaces, $\phi^{\vee} : Q^{\vee} \to
P^{\vee}$, is also $\kts$-linear.  These constructions are
functorial in $P$ and $Q$, and thus $P \mapsto P^{\vee} : \cP \to
\cP$ defines a contravariant $\FF_q(t)$-linear functor.

\subsubsection{Dual of the Carlitz motive}
Using the definition of the Carlitz motive in
\S\ref{SSS:CarlitzMotive}, we see that $C^{\vee}$ is isomorphic to
$\ok(t)$ as a $\ok(t)$-vector space and that
\[
  \bsigma f = \frac{1}{t-\theta} \cdot f^{(-1)}, \quad
  \textnormal{$f \in C^{\vee}$ ($= \ok(t)$).}
\]
Furthermore, we see that $C^{\vee} \otimes C \cong \one$ and that
$C$ is an invertible object in $\cP$.  Thus the functor
\[
  P \mapsto P \otimes C : \cP \to \cP
\]
is an equivalence of categories.  We define for $n \in \ZZ$,
\[
C(n) := \begin{cases}
  C^{\otimes n} & \textnormal{if $n > 0$,} \\
  \one & \textnormal{if $n=0$,} \\
  (C^{\vee})^{\otimes -n} & \textnormal{if $n < 0$.}
\end{cases}
\]

\subsubsection{Rigid abelian tensor category}
In the language of \cite[\S II.1]{delmil82}, it is easily shown that
the category of pre-$t$-motives is an abelian $\FF_q(t)$-linear
tensor category.  We omit the details, but we observe that
\begin{itemize}
\item each $\Hom_{\cP}(P,Q)$ is a finite dimensional vector space
over $\FF_q(t)$;
\item $\otimes$ is compatibly associative and commutative;
\item $\otimes$ is $\FF_q(t)$-bilinear;
\item $\one$ is an identity object with respect to tensor products.
\end{itemize}
Furthermore, it is straightforward to check that
\begin{itemize}
\item the pre-$t$-motive $\Hom(P,Q)$ defines an internal Hom in
$\cP$ that is compatible with tensor products;
\item for each pre-$t$-motive $P$, there is a natural isomorphism $P
\cong P^{\vee\vee}$.
\end{itemize}
Therefore, $\cP$ is also rigid.  We record this information in the
following theorem.

\begin{theorem} \label{T:PisRigidABTensor}
  The category $\cP$ of pre-$t$-motives is a rigid abelian
  $\FF_q(t)$-linear tensor category.
\end{theorem}

\subsection{Rigid analytic triviality} \label{SS:RAT}

\subsubsection{The category $\cR$} \label{SSS:RATDef}
Let $P$ be a pre-$t$-motive.  We set
\[
  \dP := \LL \otimes_{\ok(t)} P,
\]
and give $\dP$ a left $\kts$-module structure by setting
\[
  \bsigma(f \otimes m) := f^{(-1)} \otimes \bsigma m.
\]
Let
\[
  P^{\rB} := (\dP)^{\bsigma} = \{ \mu \in \dP \mid \bsigma \mu = \mu \}.
\]
Then $P^{\rB}$ is an $\FF_q(t)$-vector space, and $P \mapsto
P^{\rB}$ is a covariant functor from $\cP$ to the category of
$\FF_q(t)$-vector spaces. (The ``$\rB$'' in $P^{\rB}$ stands for
``Betti.'')  It is straightforward to check that $P \mapsto P^{\rB}$
is left exact.

We say that $P$ is \emph{rigid analytically trivial} if the natural
map
\[
  \LL \otimes_{\FF_q(t)} P^{\rB} \to \dP
\]
is an isomorphism.  If $P \cong Q$ as pre-$t$-motives and $P$ is
rigid analytically trivial, then so is $Q$.  We let $\cR$ denote the
strictly full subcategory of $\cP$ whose objects are the rigid
analytically trivial pre-$t$-motives.  Clearly the zero object is
rigid analytically trivial, and so $\cR$ is non-empty.  We shall see
momentarily that $\one$ and $C$ are also rigid analytically trivial.

\begin{lemma} \label{L:ftwist}
We have $\LL^\sigma = \FF_q(t)$.
\end{lemma}

\begin{proof}
By definition, for $f \in \LL^{\sigma}$ we have $f^{(-1)} = f$, and
so by \eqref{E:LLfactorization} the polar divisor $D$ of $f$ on the
closed unit disk in $\KK$ must also satisfy $D^{(-1)}=D$.  Therefore
$D$ is the divisor of zeros of a polynomial $c$ in $\FF_q[t]$.  Then
$cf \in \TT$, and $(cf)^{(-1)} = cf$, from which we have $cf \in \TT
\cap \power{\FF_q}{t} = \FF_q[t]$.
\end{proof}

\begin{proposition} \label{P:oneisRAT}
The pre-$t$-motive $\one$ is rigid analytically trivial.
\end{proposition}

\begin{proof}
It is clear that $\one^\dagger = \LL$ with $\bsigma f := f^{(-1)}$
for $f \in \LL$.  Therefore, by Lemma~\ref{L:ftwist}, $\one^{\rB} =
\LL^\sigma = \FF_q(t)$. Thus $\LL \otimes_{\FF_q(t)} \one^{\rB}
\cong \one^\dagger$.
\end{proof}

\subsubsection{The power series $\Omega$}
Consider the power series
\[
  \Omega = \Omega(t) := \zeta_{\theta}^{-q} \prod_{i=1}^\infty
  \Bigl( 1 - t/\theta^{(i)} \Bigr) \in
\power{k_\infty(\zeta_{\theta})}{t} \subseteq \power{\KK}{t}.
\]
It is not difficult to show that $\Omega(t)$ has an infinite radius
of convergence, and so $\Omega \in \EE \subseteq \TT$. Since
$\Omega$ has infinitely many zeros in $\KK$, it follows that $\Omega
\notin \overline{\KK(t)}$.  Since $\Omega$ has no zeros inside the
unit disk, it follows that $\Omega \in \TT^{\times}$. It also
satisfies the functional equation
\[
  \Omega^{(-1)} = (t-\theta) \Omega.
\]
The number
\[
\tpi = -\frac{1}{\Omega(\theta)} = \theta\zeta_{\theta}
\prod_{i=1}^\infty \Bigl( 1 - \theta^{1-q^i} \Bigr)^{-1} \in
k_\infty(\zeta_{\theta})
\]
is the \emph{Carlitz period}, which figures prominently in our
transcendence considerations later on (see also
\cite[Cor.~5.2.8]{andth90}, \cite[\S 3.2]{goss:FF}, \cite[\S
2.5]{thakur}).

\begin{lemma} \label{L:OmegaFE}
Suppose $f \in \LL$ satisfies $(t-\theta)^n f^{(-1)}= f$ for some $n
\in \ZZ$.  Then $f = c/\Omega^n$ for some $c \in \FF_q(t)$.
\end{lemma}

\begin{proof}
  Let $c = f\Omega^n$.  Then $c$ satisfies $c^{(-1)} = c$, and so by
  Lemma~\ref{L:ftwist}, $c \in \FF_q(t)$.
\end{proof}

\begin{proposition} \label{P:CarlitzRAT}
The Carlitz motive $C$ is rigid analytically trivial.
\end{proposition}

\begin{proof}
  We see that $C^\dagger = \LL$ with $\bsigma f = (t-\theta) f^{(-1)}$
  for $f \in \LL$.  Therefore, by Lemma~\ref{L:OmegaFE},
\[
  C^{\rB} = \{ f \in \LL \mid (t-\theta) f^{(-1)} = f \}
  = \frac{1}{\Omega} \cdot \FF_q(t).
\]
Therefore $\LL \otimes_{\FF_q(t)} C^{\rB} \cong C^\dagger$.
\end{proof}

\begin{lemma} \label{L:Bindependence}
  Let $P$ be a pre-$t$-motive, and let $\mu_1, \dots, \mu_m \in P^{\rB}$.
  If $\mu_1, \dots, \mu_m$ are linearly independent over $\FF_q(t)$, then
  they are linearly independent over $\LL$ in~$\dP$.
\end{lemma}

\begin{proof}
Suppose that $m \geq 2$ is minimal such that $\mu_1, \dots, \mu_m$
are linearly independent over $\FF_q(t)$ but that $\sum_{i=1}^m f_i
\mu_i = 0$, with $f_i \in \LL,\ f_1 = 1$. Now,
\[
  \bsigma \sum_{i=1}^m f_i \mu_i = \sum_{i=1}^m f_i^{(-1)} \mu_i = 0.
\]
Therefore, $\sum_{i=2}^m (f_i - f_i^{(-1)})\mu_i = 0$.  By the
minimality of $m$ and Lemma~\ref{L:ftwist}, each $f_i$ is in
$\FF_q(t)$. However, this violates the $\FF_q(t)$-linear
independence of $\mu_1, \dots, \mu_m$.
\end{proof}

\begin{proposition} \label{P:RATinequality}
If $P$ is a pre-$t$-motive, then $\dim_{\FF_q(t)} P^{\rB} \leq
\dim_{\ok(t)} P$.  Equality holds if and only if $P$ is rigid
analytically trivial.
\end{proposition}

\begin{proof}
{From} Lemma~\ref{L:Bindependence}, the map $\LL \otimes_{\FF_q(t)}
P^{\rB} \to \dP$ is injective.  The inequality in the statement of
the proposition follows from the equality $\dim_{\ok(t)} P =
\dim_{\LL} \dP$.  By the definition of rigid analytic triviality,
equality holds if and only if the map above is also surjective.
\end{proof}

\begin{proposition} \label{P:RATrivialization}
Suppose that $P$ is a pre-$t$-motive and that $\Phi$ represents
multiplication by $\bsigma$ on $P$ with respect to the basis $\bp$
of $P$.
\begin{enumerate}
\item[(a)] $P$ is rigid analytically trivial if and only if there is a
matrix $\Psi \in \GL_r(\LL)$ satisfying
\[
  \sigma(\Psi) = \Psi^{(-1)} = \Phi\Psi.
\]
Such a matrix $\Psi$ is called a \emph{rigid analytic
trivialization of $\Phi$} (cf.~\cite[Thm.~5]{and86},
\cite[Lem.~4.4.13]{abp04}).
\item[(b)] If $\Psi$ is a rigid analytic trivialization of $\Phi$, then the
entries of $\Psi^{-1}\bp$ form an $\FF_q(t)$-basis for $P^{\rB}$.
\item[(c)] If $P$ is rigid analytically trivial,
$\Phi \in \Mat_r(\ok[t])$, and $\det(\Phi) = d(t-\theta)^s$ for some
$s \geq 0$ and $d \in \ok^{\times}$, then there is a rigid analytic
trivialization $\Psi$ of $\Phi$ with $\Psi \in \GL_r(\TT)$.
\end{enumerate}
\end{proposition}

\begin{proof}
The proofs of parts (a) and (b) are essentially the same as the
proof of \cite[Lem.~4.4.13]{abp04} with minor modifications.  We
provide a sketch for completeness.  ((a) $\Leftarrow$; (b)):
Certainly if we have such a $\Psi$, then the entries of
$\Psi^{-1}\bp$ are both an $\LL$-basis of $\dP$ and also an
$\FF_q(t)$-linearly independent set in $P^{\rB}$.  By
Proposition~\ref{P:RATinequality}, the entries of $\Psi^{-1}\bp$
must be an $\FF_q(t)$-basis of $P^{\rB}$, and thus $P$ is rigid
analytically trivial.  ((a) $\Rightarrow$): On the other hand, if
$P$ is rigid analytically trivial, then there is a matrix $\Theta
\in \GL_r(\LL)$ so that the entries of $\Theta \bp$ are both an
$\LL$-basis of $\dP$ and an $\FF_q(t)$-basis of $P^\rB$. Setting
$\Psi := \Theta^{-1}$ gives the desired matrix.

For part (c), we first let $\sP$ be the $\ok[t]$-span of the
entries of $\bp$, and set
\[
  \sP^\dagger := \TT \otimes_{\ok[t]} \sP, \quad \sP^{\rB} := \{
  \mu \in \sP^\dagger \mid \bsigma \mu = \mu \}.
\]
For $\mu \in P^\rB$, write $\mu = \sum f_i p_i = \bff\cdot \bp$ with
$\bff \in \Mat_{1\times r}(\LL)$.  We claim that for some $c \in
\FF_q[t]$, we have $c\mu \in \sP^\rB$.  Let $\den(\bff) \in \KK[t]$
denote the monic least common multiple of the denominators of
$\bff$, which is well-defined by \eqref{E:LLfactorization}. Then
since $\bsigma \mu = \mu$, we have
\[
  \bff \cdot \bp = \bsigma(\bff \cdot \bp) = \bff^{(-1)} \cdot \Phi \bp.
\]
Therefore, $\den(\bff) = \den(\bff^{(-1)} \cdot \Phi)$.  But $\Phi
\in \Mat_r(\ok[t])$, so $\den(\bff^{(-1)} \cdot \Phi)$ divides
$\den(\bff^{(-1)})$.  Degree considerations force $\den(\bff) =
\den(\bff^{(-1)})$.  Therefore take $c = \den(\bff) \in \FF_q[t]$.
This proves the claim, and moreover we have shown that
\[
  P^{\rB} \cong \FF_q(t) \otimes_{\FF_q[t]} \sP^{\rB}.
\]
Furthermore, it follows that as $\LL$-vector spaces, $\dP \cong \LL
\otimes_{\TT} \sP^\dagger \cong \LL \otimes_{\FF_q[t]} \sP^{\rB}$.
Let $\bnu = [\nu_1, \dots, \nu_r ]^{\tr}$ be an $\FF_q[t]$-basis for
$\sP^{\rB}$. Then for some $\Theta \in \GL_r(\LL) \cap \Mat_r(\TT)$,
we have $\bnu = \Theta \bp$.  Since $\bsigma \bnu = \bnu$, it
follows that $\Theta^{(-1)} \Phi = \Theta$.  By our initial
hypotheses, $d(t-\theta)^s \det(\Theta)^{(-1)} = \det(\Theta)$.
Choose $b \in \ok^{\times}$ so that $d=b^{(-1)}/b$. Then from
Lemmas~\ref{L:ftwist} and~\ref{L:OmegaFE} (and the fact that $\Theta
\in \Mat_r(\TT)$), we see that
\[
  b\det(\Theta) = \frac{\gamma}{\Omega^s}, \quad \gamma \in
  \FF_q[t].
\]
We claim that $\gamma \in \FF_q^{\times}$.  If not, then
$\det(\Theta) \equiv 0 \pmod{\gamma}$ in $\TT$, and so there is a
$\bff = [f_1, \dots, f_r] \in \Mat_{1\times r}(\TT)$ so that
\[
  \bff \cdot \Theta \equiv 0 \pmod{\gamma}.
\]
Since $\TT/\gamma\TT \cong \KK[t]/\gamma\KK[t]$, without loss of
generality we can assume that each $f_i$ is a polynomial in $\KK[t]$
of degree strictly less than the degree of $\gamma$, that
$\absT{f_i}\leq 1$ for all $i$, and that at least one $f_i$
satisfies $\absT{f_i}=1$.  Now define a norm
$\absT{\,\cdot\,}_{\dagger}$ on $\sP^\dagger$ by
\[
  \biggl\| \sum h_i p_i \biggr\|_{\dagger} := \sup \absT{h_i},
\quad h_1, \dots, h_r \in \TT.
\]
Then $\absT{\,\cdot\,}_\dagger$ defines a complete ultrametic norm
on $\sP^\dagger$ that satisfies
\[
  \absT{h \mu}_{\dagger} = \absT{h} \cdot \absT{\mu}_{\dagger},
  \quad h \in \TT,\ \mu \in \sP^\dagger.
\]
Consider
\[
  \bff \cdot \Theta^{(-1)} \Phi = \bff\cdot \Theta \equiv 0 \pmod{\gamma}.
\]
Since $\gamma$ is relatively prime to $\det(\Phi)$, it follows that
$\Phi$ is invertible modulo $\gamma$, and so
\[
  \bff \cdot \Theta^{(-1)} \equiv 0 \pmod{\gamma}.
\]
Repeating this argument we find that $\bff\cdot \Theta^{(-n)} \equiv
0 \pmod{\gamma}$ for all $n \geq 0$.  Now, by choice of $\bff$,
\[
  \frac{1}{\gamma} \bff \cdot \bnu = \frac{1}{\gamma} \bff \cdot \Theta\cdot
   \bp \in \sP^\dagger,
\]
and for each $n$, the above congruences for $\bff \cdot
\Theta^{(-n)}$ imply that $\gamma^{-1} \bff^{(n)} \cdot \bnu =
\gamma^{-1} \bff^{(n)} \cdot \Theta \cdot \bp \in \sP^\dagger$. Now
by Lemma~\ref{L:LLTwist}, there is an $m
>0$ so that with respect to the $\absT{\,\cdot\,}_\dagger$ metric,
\[
  \lim_{n \to \infty} \frac{1}{\gamma} \sum f_i^{(mn)} \nu_i =
  \frac{1}{\gamma} \sum c_i \nu_i \in \sP^\dagger,
\]
where $c_i \in \overline{\FF}_q[t]$ and at least one $c_i \neq 0$,
say $c_a \neq 0$.  Now for some $l \geq 1$, we have every $c_i \in
\FF_{q^l}[t]$.  Since the trace map $\FF_{q^l} \to \FF_q$ is not
trivial, by dividing each $c_i$ by a fixed element in
$\FF_{q^l}^\times$, we can assume that $c_a + c_a^{(-1)} + \dots +
c_a^{(1-l)} \neq 0$.  Therefore,
\[
  \sum_{j=0}^{l-1} \bsigma^{j} \biggl( \frac{1}{\gamma}
   \sum_{i=1}^r c_i \nu_i \biggr) = \frac{1}{\gamma} \sum_{i=1}^r
   \biggl( \sum_{j=0}^{l-1} c_i^{(-j)} \biggr) \nu_i
  \in \sP^\dagger.
\]
Thus we obtain $\mu := \gamma^{-1} \sum d_i \nu_i \in \sP^\dagger$,
$d_i \in \FF_q[t]$, $d_a \neq 0$. Easily we see that $\mu \in
\sP^{\rB}$ and $\mu \neq 0$.  Since $\deg f_i < \deg \gamma$ for
each $i$, we have $\deg d_i < \deg \gamma$ for each $i$. In
particular, $\gamma$ does not divide $d_a$. Thus $\mu \in \sP^{\rB}$
but $\mu$ is not in the $\FF_q[t]$-span of $\bnu$, which contradicts
that $\bnu$ is an $\FF_q[t]$-basis of $\sP^{\rB}$. Therefore, it
follows that $\gamma \in \FF_q^{\times}$, and since $\Omega \in
\TT^{\times}$, we have $\det(\Theta) \in \TT^{\times}$. Taking $\Psi
= \Theta^{-1}$ provides the desired rigid analytic trivialization.
\end{proof}

\subsubsection{Remark} \label{SSS:RATiffRAT}
It is worth noting that multiplication by
$\Theta$ induces the isomorphism of $\LL$-vector spaces,
\[
  \LL \otimes_{\TT} \sP^{\dagger} \cong \LL \otimes_{\TT} (\TT
  \otimes_{\FF_q[t]} \sP^{\rB}).
\]
Since $\Theta \in \GL_r(\TT)$, this then implies $\sP^{\dagger}
\cong \TT \otimes_{\FF_q[t]} \sP^{\rB}$ as $\TT$-modules.

\begin{proposition} \label{P:RATKerQuot}
Let
\[
  0 \to P \to Q \to R \to 0
\]
be an exact sequence of pre-$t$-motives.
\begin{enumerate}
\item[(a)] If $Q$ is rigid analytically trivial, then both $P$ and $R$
  are rigid analytically trivial.
\item[(b)] If $P$, $Q$, and $R$ are rigid analytically trivial, then
  the sequence
\[
0 \to P^\rB \to Q^\rB \to R^\rB \to 0
\]
is an exact sequence of $\FF_q(t)$-vector spaces.
\end{enumerate}
\end{proposition}

\begin{proof}
  The sequence $0 \to P^{\rB} \to Q^{\rB} \to R^{\rB}$ is exact.  Now
  suppose that $Q$ is rigid analytically trivial.  Let $\kappa : \LL
  \otimes_{\FF_q(t)} Q^{\rB} \to \LL \otimes_{\FF_q(t)} R^{\rB}$ be the
  natural map.  Then we have a commutative diagram with exact rows,
\[
\xymatrix{
0 \ar[r] & \LL \otimes_{\FF_q(t)} P^{\rB} \ar[r] \ar[d] &
  \LL \otimes_{\FF_q(t)} Q^{\rB} \ar[r] \ar[d]^{\wr} &
  \im(\kappa) \ar[r] \ar[d] & 0 \\
0 \ar[r] & \dP \ar[r] & \dQ \ar[r] & R^\dagger \ar[r]& 0,
}
\]
where the central vertical map is an isomorphism by hypothesis,
and the other two are injective by Lemma~\ref{L:Bindependence}.
The injectivity of all three maps then implies that each is an
isomorphism. Thus we see immediately that $P$ is rigid
analytically trivial.  Also we see that
\[
\dim_{\FF_q(t)} R^{\rB} = \dim_{\LL} \LL \otimes_{\FF_q(t)}
R^{\rB}
\geq \dim_{\LL} \im(\kappa) = \dim_{\LL} R^{\dagger} =
\dim_{\ok(t)} R,
\]
which by Proposition~\ref{P:RATinequality} must be a string of
equalities.  Therefore $R$ is rigid analytically trivial, which
completes part (a).

Now suppose that $P$, $Q$, and $R$ are all rigid analytically trivial.  Then
\[
\dim_{\FF_q(t)} Q^\rB = \dim_{\ok(t)} Q = \dim_{\ok(t)} P +
\dim_{\ok(t)} R = \dim_{\FF_q(t)} P^\rB + \dim_{\FF_q(t)} R^\rB,
\]
which proves part (b).
\end{proof}

\subsubsection{Remark} \label{SSS:RATAbBExact}
  In particular, it follows from Proposition~\ref{P:RATKerQuot} that
  kernels and cokernels exist in $\cR$, which implies that $\cR$ is
  an abelian $\FF_q(t)$-linear category.   We also see that
\[
  P \to P^{\rB} : \cR \to \Vector{\FF_q(t)},
\]
where $\Vector{\FF_q(t)}$ is the category of finite dimensional
vector spaces over $\FF_q(t)$, is an exact $\FF_q(t)$-linear
functor.

\begin{proposition} \label{P:Bfaithful}
  Let $P$ and $Q$ be rigid analytically trivial pre-$t$-motives.  Then
  the natural map
\[
  \Hom_{\cR}(P,Q) \to \Hom_{\FF_q(t)}(P^{\rB},Q^{\rB})
\]
is injective.
\end{proposition}

\begin{proof}
  Suppose $\phi : P \to Q$ is a morphism in $\Hom_{\cR}(P,Q)$.  Then
  we have an exact sequence in $\cR$,
\[
  0 \to \ker \phi \to P \stackrel{\phi}{\to} Q \to Q/\phi(P) \to 0,
\]
which leads then to an exact sequence of $\FF_q(t)$-vector spaces,
\[
  0 \to (\ker \phi)^{\rB} \to P^{\rB} \stackrel{\phi^{\rB}}{\to}
  Q^{\rB} \to (Q/\phi(P))^{\rB} \to 0.
\]
Since the dimension over $\ok(t)$ of each term in the first sequence
is the same as the dimension over $\FF_q(t)$ of the corresponding term
in the second sequence, we see that $\phi^{\rB} = 0$ if and only if
$\phi = 0$.
\end{proof}

\begin{proposition} \label{P:RATTensorDual}
  If pre-$t$-motives $P$ and $Q$ are rigid analytically trivial, then
\begin{enumerate}
\item[(a)] $P \otimes Q$ is rigid analytically trivial, and the natural map
$P^{\rB} \otimes_{\FF_q(t)} Q^{\rB} \to (P \otimes Q)^{\rB}$ is an
isomorphism of $\FF_q(t)$-vector spaces;
\item[(b)] $P^\vee$ is rigid analytically trivial, and the natural map
$\bigl(P^{\rB}\bigr)^{\vee} \to \bigl(P^{\vee}\bigr)^{\rB}$ is an
isomorphism of $\FF_q(t)$-vector spaces.
\end{enumerate}
\end{proposition}

\begin{proof}
Here we make use of Proposition~\ref{P:RATrivialization}.   We
first note that
\[
  (P \otimes Q)^{\dagger} = \LL \otimes_{\ok(t)} ( P \otimes_{\ok(t)} Q)
\cong (\LL \otimes_{\ok(t)} P) \otimes_{\LL} (\LL \otimes_{\ok(t)} Q)
= \dP \otimes_{\LL} \dQ,
\]
where the middle isomorphism is an isomorphism of $\LL$-vector
spaces that commutes with the action of $\bsigma$.  We observe that
we can choose $\ok(t)$-bases for $P$, $Q$, and $P\otimes Q$ so that
multiplication by $\bsigma$ is represented by matrices $\Phi_P$,
$\Phi_Q$, and $\Phi_{P\otimes Q}$ satisfying
\[
  \Phi_{P \otimes Q} = \Phi_P \otimes \Phi_Q.
\]
By Proposition~\ref{P:RATrivialization}, we can choose $\Psi_P$,
$\Psi_Q \in \GL_r(\LL)$ that are rigid analytic trivializations of
$\Phi_P$ and $\Phi_Q$.  Then we note that $\Psi_{P\otimes Q} :=
\Psi_P \otimes \Psi_Q$ is a rigid analytic trivialization of
$\Phi_{P \otimes Q}$.  Now  note that $\Phi_{P^\vee} :=
(\Phi_P^{-1})^{\tr}$ represents multiplication by $\bsigma$ with
respect to the dual basis and that $\Psi_P^{\vee} :=
(\Psi_P^{-1})^{\tr}$ is a rigid analytic trivialization.  The second
parts of (a) and (b) are straightforward.
\end{proof}

\begin{theorem}
The category $\cR$ of rigid analytically trivial pre-$t$-motives is
a neutral Tannakian category over $\FF_q(t)$ with fiber functor $P
\mapsto P^{\rB} : \cR \to \Vector{\FF_q(t)}$.
\end{theorem}

\begin{proof}
We have seen that
\begin{itemize}
\item $\one$ is in $\cR$ (Proposition~\ref{P:oneisRAT});
\item $\cR$ is an abelian category (Proposition~\ref{P:RATKerQuot} and
  \S\ref{SSS:RATAbBExact});
\item $\cR$ is closed under tensor products and duals
  (Proposition~\ref{P:RATTensorDual}).
\end{itemize}
Thus $\cR$ is a rigid abelian $\FF_q(t)$-linear tensor subcategory of
$\cP$ (see \cite[Defs.~II.1.14-15]{delmil82}).  We have also shown
that
\begin{itemize}
\item $\End_{\cR}(\one) = \FF_q(t)$ (Lemma~\ref{L:End1});
\item For each $P$ in $\cR$, the $\FF_q(t)$-vector space $P^\rB$
is finite dimensional (Proposition~\ref{P:RATinequality});
\item $P \mapsto P^{\rB}$ is $\FF_q(t)$-linear and exact
  (Proposition~\ref{P:RATKerQuot} and \S\ref{SSS:RATAbBExact});
\item $P \mapsto P^{\rB}$ is faithful (Proposition~\ref{P:Bfaithful});
\item $P \mapsto P^{\rB}$ is a tensor functor
  (Proposition~\ref{P:RATTensorDual}).
\end{itemize}
Thus $\cR$ is a neutral Tannakian category over $\FF_q(t)$ with fiber
functor $P \mapsto P^{\rB}$ (see \cite[Def.~II.2.19]{delmil82}).
\end{proof}

\subsection{Anderson $t$-motives} \label{SS:AndersontMotives}

Here we recall the definitions and essential properties of ``dual
$t$-motives'' from \cite{abp04}.  So as not to confuse these
objects with the duals of $t$-motives to be used later on, we call
these objects \emph{Anderson $t$-motives}, since they are simply
the dual notion of the objects studied in \cite{and86}.

\subsubsection{Definition} \label{SSS:AndersontMotiveDefinition}
An \emph{Anderson $t$-motive} $\sM$ is a left $\ktsp$-module such
that
\begin{itemize}
\item $\sM$ is free and finitely generated over $\ok[t]$;
\item $\sM$ is free and finitely generated over $\ok[\bsigma]$;
\item $(t-\theta)^n \sM \subseteq \bsigma \sM$ for all $n \gg 0$.
\end{itemize}
A morphism of Anderson $t$-motives is a left $\ktsp$-module
homomorphism.  In this way Anderson $t$-motives form a category.

As in \S\ref{SSS:RepsPretmot}, if $\sm \in \Mat_{r
  \times 1}(\sM)$ is a $\ok[t]$-module basis for $\sM$, then there is
a matrix $\Phi = \Phi_{\sm} \in \Mat_{r \times 1}(\ok[t])$ so that
\[
  \bsigma \sm = \Phi \sm.
\]
Since a power of $t-\theta$ annihilates $\sM/\bsigma \sM$, we have
\[
  \det \Phi = c(t-\theta)^s
\]
for some $c \in \ok^{\times}$, where $s$ is the rank of $\sM$ as a
$\ok[\bsigma]$-module.

\subsubsection{Anderson $t$-motives to pre-$t$-motives}
Given an Anderson $t$-motive $\sM$ we obtain a pre-$t$-motive $M$
by setting
\[
  M := \ok(t) \otimes_{\ok[t]} \sM
\]
and defining
\[
  \bsigma (f \otimes m) := f^{(-1)} \otimes \bsigma m.
\]
It is straightforward to check that $M$ is a left $\kts$-module, and
it is of course finite dimensional as a $\ok(t)$-vector space.
Moreover, $\sM \mapsto M$ is a functor from the category of Anderson
$t$-motives to the category of pre-$t$-motives.

\subsubsection{The Carlitz motive}
\label{SSS:AndersonCarlitzMotive} Let $\sC$ be the Anderson
$t$-motive whose underlying $\ok[t]$-module is $\ok[t]$ itself.
Then the action of $\bsigma$ on $\sC$ is defined by
\[
  \bsigma(f) = (t-\theta)f^{(-1)}, \quad f \in \sC.
\]
It is not difficult to check that $\sC$ is an Anderson $t$-motive,
and that its image in $\cP$ is the Carlitz motive.  For any $n
\geq 1$, we also have the $n$-th tensor power of $\sC$,
\[
  \sC(n) := \sC \otimes_{\ok[t]} \dots \otimes_{\ok[t]}
  \sC,
\]
with diagonal $\bsigma$-action.  It is an Anderson $t$-motive sent
to $C(n)$ in $\cP$.

\subsubsection{The Carlitz module} \label{SSS:CarlitzModule}
The \emph{Carlitz module $\fC$ over $\ok$} is defined to be the
$\FF_q$-algebra $\ok$ together with an $\FF_q[t]$-module structure
defined by
\[
  \fC_t(x) := \theta x + x^q, \quad x \in \ok.
\]
That is the $\FF_q$-algebra homomorphism $a \mapsto \fC_a : \FF_q[t]
\to \ok[\bsigma^{-1}]$ defined by $t \mapsto \theta + \bsigma^{-1}$
induces an $\FF_q[t]$-module structure on $\ok$.  See
\cite[Ch.~3]{goss:FF} or \cite[\S 2.5]{thakur} for more details.  To
see the relationship with the Carlitz motive, we note that there is
an isomorphism
\[
  \fC(\ok) \cong \frac{\sC}{(\bsigma-1)\sC}
\]
of $\FF_q[t]$-modules.  Indeed if $x \in \ok$, then
\[
tx = \theta x + (t-\theta)x = \theta x + \bsigma(x^q) = \theta x +
x^q + (\bsigma -1)x^q.
\]
Similarly $ax \equiv \fC_a(x) \pmod{\bsigma-1}$ for all $a \in
\FF_q[t]$.  It is a simple matter to check that there is a natural
isomorphism of $\FF_q$-vector spaces $\sC/(\bsigma-1)\sC \cong \ok$.
Thus $\sC/(\bsigma-1)\sC$ presents the Carlitz module directly.

\begin{proposition} \label{P:AtoPFullyFaithful}
For Anderson $t$-motives $\sM$ and $\sN$, the natural map
\[
  \Hom_{\ktsp}(\sM,\sN) \otimes_{\FF_q[t]} \FF_q(t) \to \Hom_{\cP}(M,N)
\]
is an isomorphism of $\FF_q(t)$-vector spaces.
\end{proposition}

\begin{proof}
Let $\Theta$ denote the map in question.  It is clearly
$\FF_q(t)$-linear.  To see that it is injective, we first observe
that if $\alpha \in \Hom_{\ktsp}(\sM,\sN) \otimes_{\FF_q[t]}
\FF_q(t)$ then $\alpha = \uphi \otimes \frac{1}{v}$, for some $\uphi
\in \Hom_{\ktsp}(\sM,\sN)$ and $v \in \FF_q[t]$, $v \neq 0$. Then $v
\Theta(\alpha) = \Theta(v\alpha) = \Theta(\uphi \otimes 1) =: \phi$.
But
\begin{align*}
\uphi &\in \Hom_{\ktsp}(\sM,\sN) \subseteq
\Hom_{\ok[t]}(\sM,\sN), \\
\phi &\in \Hom_{\kts}(M,N) \subseteq \Hom_{\ok(t)}(M,N),
\end{align*}
and so $\phi = 0$ if and only if $\uphi = 0$.  Thus
$\Theta(\alpha) = 0$ if and only if $\alpha = 0$.

For surjectivity, suppose that $\phi \in \Hom_{\cP}(M,N)$. Fix
$\ok[t]$-bases $\sm$ and $\sn$ for $\sM$ and $\sN$ respectively, and
extend these to bases $\bm$ and $\bn$ of $M$ and $N$. Then the map
$\phi : M \to N$ is represented by a matrix $F \in \Mat_{r \times
s}(\ok(t))$ so that $F^{(-1)} \Phi_{\bn} = \Phi_{\bm} F$ as in
\S\ref{SSS:RepsPretmot}.  By choice of $\sm$ and $\sn$, $\Phi_{\bm}$
and $\Phi_{\bn}$ have entries in $\ok[t]$, and it suffices to show
that $F$ has entries with denominators in $\FF_q[t]$.

For a matrix $B$ with entries in $\ok(t)$, let $\den(B) \in
\ok[t]$ be the monic least common multiple of the denominators of
the entries of $B$.  Since $\det(\Phi_{\bn}) = c(t-\theta)^s$ for
some $s \geq 0$ and $c \in \ok^{\times}$, we see that
\[
  \den(F) (t-\theta)^s \cdot F^{(-1)} = \den(F) (t-\theta)^s \cdot
  \Phi_{\bm}\,F\, \Phi_{\bn}^{-1} \in \Mat_{r\times s}(\ok[t]).
\]
Therefore, $\den(F^{(-1)})$ divides $\den(F)(t-\theta)^s$. However,
$\den(F^{(-1)}) = \den(F)^{(-1)}$ and so $\deg(\den(F^{(-1)})) =
\deg(\den(F))$.  Thus, it suffices to show that $\den(F^{(-1)})$ is
relatively prime to $t-\theta$, since then $\den(F)^{(-1)} =
\den(F)$ whence all of the denominators of $F$ are in $\FF_q[t]$.

Suppose that $t-\theta$ divides $\den(F^{(-1)})$, and so
$t-\theta^q$ divides $\den(F)$.  Then $t-\theta^q$ divides
$\den(\Phi_{\bm}F)$, because otherwise $t-\theta^q$ would divide
$\det(\Phi_{\bm})$ which is a power of $t-\theta$. Likewise,
$t-\theta^q$ divides $\den(\Phi_{\bm}F\Phi_{\bn}^{-1}) =
\den(F^{(-1)})$.  By repeating the same argument we see that
$\den(F^{(-1)})$ is divisible by each of
\[
  t-\theta, t - \theta^q, t - \theta^{q^2}, \ldots
\]
contradicting that $\den(F^{(-1)}) \in \ok[t]$.
\end{proof}

\subsubsection{Rigid analytic triviality}
Similar to \S\ref{SSS:RATDef}, if $\sM$ is an Anderson $t$-motive,
then we set
\[
  \sM^\dagger := \TT \otimes_{\ok[t]} \sM.
\]
We provide $\sM^\dagger$ with a $\ktsp$-module structure by setting
$\bsigma(f \otimes m) = f^{(-1)} \otimes \bsigma m$, and we set
\[
  \sM^{\rB} := (\sM^\dagger)^{\bsigma} =
  \{ \mu \in \sM^\dagger \mid \bsigma \mu = \mu \}.
\]
We say that $\sM$ is \emph{rigid analytically trivial} if the
natural map $\TT \otimes_{\FF_q[t]} \sM^{\rB} \to \sM^\dagger$ is an
isomorphism.  The following proposition is a companion to
Proposition~\ref{P:RATrivialization}.

\begin{proposition} \label{P:ARATrivialization}
Let $\sM$ be an Anderson $t$-motive, and let $M$ be its
corresponding pre-$t$-motive.  Suppose $\sm \in \Mat_{r\times
1}(\sM)$ is a $\ok[t]$-basis for $\sM$, and let $\Phi \in
\Mat_{r\times 1}(\ok[t])$ represent multiplication by $\bsigma$ on
$\sM$ with respect to $\sm$.
\begin{enumerate}
\item[(a)] $\sM$ is rigid analytically trivial if and only if
it admits a rigid analytic trivialization $\Psi$ with $\Psi \in
\GL_r(\TT)$.
\item[(b)] If $\Psi \in \GL_r(\TT)$ is a rigid analytic
trivialization of $\Phi$, then the entries of $\Psi^{-1} \sm$ form
an $\FF_q[t]$-basis of $\sM^{\rB}$.
\item[(c)] $\sM$ is rigid analytically trivial if and only if
$M$ is rigid analytically trivial.
\end{enumerate}
\end{proposition}

\begin{proof}
The proofs of parts (a) and (b) are in \cite[Lem.~4.4.13]{abp04}
and follow the same lines as their counterparts in
Proposition~\ref{P:RATrivialization}.  Part (c) is then a
consequence of Proposition~\ref{P:RATrivialization}(c).
\end{proof}

\subsubsection{Definition}
We define the \emph{category $\cA^I$ of Anderson $t$-motives up to
isogeny} as follows:
\begin{itemize}
\item Objects of $\cA^I$: Anderson $t$-motives;
\item Morphisms of $\cA^I$: For
Anderson $t$-motives $\sM$ and $\sN$,
\[
\Hom_{\cA^I}(\sM,\sN) := \Hom_{\ktsp}(\sM,\sN) \otimes_{\FF_q[t]}
\FF_q(t).
\]
\end{itemize}
We also define the full subcategory $\cA\cR^I$ of rigid
analytically trivial Anderson $t$-motives up to isogeny by
restriction.  We sum up the results of this section in the
following theorem.

\begin{theorem} \label{T:AndersonIsogenyFullyFaithful}
Let $\cP$ be the category of pre-$t$-motives, and let $\cR$ be the
category of rigid analytically trivial pre-$t$-motives.
\begin{enumerate}
\item[(a)] The functor $\sM \mapsto M: \cA^I \to \cP$ is fully faithful.
\item[(b)] The functor $\sM \mapsto M: \cA\cR^I \to \cR$ is fully
faithful.
\end{enumerate}
\end{theorem}

\begin{proof}
Part (a) is simply a restatement of
Proposition~\ref{P:AtoPFullyFaithful}.  That the functor in part (b)
is well-defined follows from
Proposition~\ref{P:ARATrivialization}(c), and its full faithfulness
follows from Proposition~\ref{P:AtoPFullyFaithful}.
\end{proof}

\subsubsection{The category $\cT$}
\label{SSS:tMotiveDefinition} We define the \emph{category $\cT$ of
$t$-motives} to be the strictly full Tannakian subcategory of $\cR$
generated by the essential image of the functor
\[
  \sM \mapsto M : \cA\cR^I \to \cR.
\]
The category of $t$-motives can further be described as follows:
\begin{itemize}
\item Objects of $\cT$: rigid analytically trivial pre-$t$-motives
that can be constructed from Anderson $t$-motives using direct
sums, subquotients, tensor products, duals, and internal Hom's.
\item Morphisms of $\cT$: morphisms of left $\kts$-modules.
\end{itemize}
It is worth noting that Proposition~\ref{P:ARATrivialization}(c)
says that the category of $t$-motives is the strictly full
Tannakian subcategory of $\cR$ generated by the intersection in
$\cP$ of $\cR$ and the image of \emph{all} Anderson $t$-motives.

\subsection{Galois groups of $t$-motives}
\label{SS:GaloisGroupstMotives}

Having defined a Tannakian category of $t$-motives, it is now
possible to assign to each $t$-motive a linear algebraic group
over $\FF_q(t)$, which we call the \emph{Galois group} of the
$t$-motive.  For essential facts about Tannakian categories and
their associated groups, we refer to \cite{breen94},
\cite{delmil82}, \cite[App.~B]{vdpsing03}.

\subsubsection{Fiber functors}
The functor
\[
\begin{array}{rcl}
\omega : \cT &\to& \Vector{\FF_q(t)} \\
M &\mapsto& M^{\rB}
\end{array}
\]
is the fiber functor of $\cT$.  For any commutative
$\FF_q(t)$-algebra $R$, we let $\omega^{(R)} : \cT \to \Module{R}$
be the extension of $\omega$ defined by
\[
  \omega^{(R)}(M) := R \otimes_{\FF_q(t)} M^{\rB},
\]
where $\Module{R}$ is the category of finitely generated left
$R$-modules.  Now fix a $t$-motive $M$. We let $\cT_M$ be the
strictly full Tannakian subcategory of $\cT$ generated by $M$. That
is, $\cT_M$ consists of all objects of $\cT$ isomorphic to
subquotients of finite direct sums of $M^{\otimes u} \otimes
(M^{\vee})^{\otimes v}$ for various $u$, $v$. The fiber functor of
$\cT_M$ is $\omega_M : \cT_M \to \Vector{\FF_q(t)}$, the restriction
of $\omega$ to $\cT_M$, and similarly we restrict $\omega_M^{(R)}$
to $\cT_M$ for an $\FF_q(t)$-algebra $R$.

\subsubsection{Galois groups} \label{SSS:tMotiveGaloisGroups}
As $\cT$ is a neutral Tannakian category over $\FF_q(t)$, there is
an affine group scheme $\Gamma_{\cT}$ over $\FF_q(t)$ so that
$\cT$ is equivalent to the category $\Rep{\Gamma_{\cT}}{\FF_q(t)}$
of finite dimensional representations of $\Gamma_{\cT}$ over
$\FF_q(t)$:
\[
  \cT \approx \Rep{\Gamma_{\cT}}{\FF_q(t)}.
\]
The group $\Gamma_{\cT}$ is defined to be the group of tensor
automorphisms of the fiber functor $\omega$; that is, if $R$ is
any $\FF_q(t)$-algebra, then
\[
  \Gamma_{\cT}(R) = \Aut^{\otimes}_{\cT}\bigl(\omega^{(R)}\bigr).
\]
Now for any $t$-motive $M$, there is a linear algebraic group
$\Gamma_M := \Gamma_{\cT_M}$ over $\FF_q(t)$ so that $\cT_M$ is
equivalent to $\Rep{\Gamma_M}{\FF_q(t)}$.  As such, for any
$\FF_q(t)$-algebra $R$, $\Gamma_M(R) =
\Aut_{\cT_M}^{\otimes}\bigl(\omega_M^{(R)}\bigr)$.  In this way we
find that we have a naturally defined faithful representation
\[
  \Gamma_M \hookrightarrow \GL(M^{\rB})
\]
over $\FF_q(t)$, which provides the basis for constructing the
equivalence of categories,
\[
  \cT_M \approx \Rep{\Gamma_M}{\FF_q(t)}.
\]
The group $\Gamma_M$ is called the \emph{Galois group of $M$}.
Furthermore, there is a surjective group homomorphism,
\[
  \Gamma_{\cT} \twoheadrightarrow \Gamma_M.
\]
If $N$ is another $t$-motive in $\cT_M$, then there is a natural
surjective homomorphism, $\Gamma_M \twoheadrightarrow \Gamma_N$. In
\S\ref{S:Galois}, we will show that $\Gamma_M$ can be calculated
using systems of $\sigma$-semilinear equations.  For now we will
calculate the Galois group of the Carlitz motive~$C$.

\begin{lemma} \label{L:HomCmCn}
For $m$, $n \in \ZZ$,
\[
  \Hom_{\cT} (C(m),C(n)) \cong
  \begin{cases}
  \FF_q(t) & \textnormal{if $m=n$,} \\
  0 & \textnormal{if $m \neq n$.}
  \end{cases}
\]
\end{lemma}

\begin{proof}
By tensoring with $C(-m)$, we see that $\Hom_{\cT}(C(m),C(n)) \cong
\Hom_{\cT}(\one,C(n-m))$.  Thus it suffices to assume that $m=0$. If
$\phi : \one \to C(n)$ is a morphism in $\cT$, then $\phi$ is
represented by some $a \in \ok(t)^\times$ such that $a =
a^{(-1)}(t-\theta)^n$.  By Lemma~\ref{L:OmegaFE}, this equation has
no non-zero solutions $a \in \ok(t)$ unless $n=0$, in which case $a$
can be anything in $\FF_q(t)$.
\end{proof}

\begin{theorem} \label{T:GaloisGroupCarlitz}
For the Carlitz motive $C$, there is an isomorphism $\Gamma_C
\cong \Gm$ over $\FF_q(t)$.
\end{theorem}

\begin{proof}
It is easy enough to check this theorem directly.  However, by
Lemma~\ref{L:HomCmCn}, $\cT_C$ is equivalent to a $\ZZ$-graded
category of vector spaces over $\ok(t)$ with a fiber functor to
$\Vector{\FF_q(t)}$, and so its Galois group is $\Gm$ over
$\FF_q(t)$ \cite[Ex.~II.2.30]{delmil82}.
\end{proof}

\section{Galois theory of systems of $\sigma$-semilinear
equations} \label{S:Galois}

In this section we demonstrate how to calculate the Galois group of
a $t$-motive as the Galois group of a system of difference equations
with respect to the automorphism $\sigma : \ok(t) \to \ok(t)$.
These systems of equations and their Galois groups are similar to
systems of linear differential equations and their Galois groups,
and one should compare our constructions with
\cite[Ch.~1]{vdpsing97}, \cite[Chs.~1--2]{vdpsing03}, which we have
used as guides, as well as \cite{andre01}, \cite{beuk92},
\cite{del90}, \cite{kolchin}, \cite{magid}, \cite{vdp99}.  For an
example of a Galois group of this type in the context of
$t$-motives, see also the proof of~\cite[Prop.~7.1]{boehar}.

Van der Put and Singer~\cite{vdpsing97} have developed the theory of
Picard-Vessiot rings for linear difference equations which is quite
useful in our context. However, their treatment generally assumes
that the field of constants is algebraically closed.  In our case
the field of constants is $\FF_q(t)$, which presents several
difficulties. On the other hand, the Picard-Vessiot rings treated in
\cite{vdpsing97} are not always domains, whereas our central
Picard-Vessiot rings \emph{are} domains by construction, which
provides several benefits for the characterization of their Galois
groups. It is worth noting that some of what is covered here is
covered by the theory of Y.~Andr\'e \cite{andre01}, but we present
everything from scratch for completeness.  We thank the referee for
making several useful suggestions for improving the clarity of this
section.

\subsection{Solutions of $\sigma$-semilinear equations}
\label{SS:sigmaSemilinearEquations}

\subsubsection{Fields of definition}
Let $F \subseteq K \subseteq L$ be fields together with an
automorphism $\sigma : L \to L$.  We say that the triple $(F,K,L)$
is \emph{$\sigma$-admissible} if
\begin{itemize}
\item $\sigma$ restricts to automorphisms of $F$ and $K$;
\item $F = F^{\sigma} = K^{\sigma} = L^{\sigma}$;
\item $L$ is a separable extension of $K$.
\end{itemize}
The primary example of $\sigma$-admissible fields that we have in
mind is
\[
  (F,K,L) = (\FF_q(t), \ok(t), \LL),
\]
with automorphism $\sigma$ defined as in \S\ref{SSS:Twisting} by
$\sigma(f) = f^{(-1)}$.  This example will be important for
applications to $t$-motives in \S\ref{SS:GammaPsiAndtMotives}. To
see that this triple is $\sigma$-admissible, we know that
$\FF_q(t)=\FF_q(t)^\sigma = \ok(t)^\sigma$ by definition and that
$\LL^\sigma=\FF_q(t)$ by Lemma~\ref{L:ftwist}.  Also, since $\LL$ is
linearly disjoint from $\ok(t^{1/p})$, it is therefore separable
over $\ok(t)$~\cite[Thm.\ 26.3]{mats2}.  Henceforth we shall assume
that a $\sigma$-admissible triple $(F,K,L)$ has been chosen.

\subsubsection{Convention}
If $\rho : S \to R$ is a homomorphism of modules or rings, and $B
\in \Mat_{r \times s}(S)$, we let $\rho(B) \in \Mat_{r\times s}(R)$
be the matrix obtained by applying $\rho$ to the entries of $B$.

\subsubsection{Definition}
Given a matrix $\Phi \in \GL_r(K)$, we consider vectors $\psi \in
\Mat_{r \times 1}(L)$ that satisfy
\[
  \sigma(\psi) = \Phi \psi.
\]
In this way, we define a \emph{system of $\sigma$-semilinear
equations}, and $\psi$ is a solution.  The set of solutions
\[
  \Sol(\Phi) := \{ \psi \in \Mat_{r \times 1}(\bL) \mid \sigma(\psi)
  = \Phi \psi \}
\]
is an $F$-vector space.

\begin{lemma} \label{L:SolIndependence}
Let $\Phi \in \GL_r(K)$.  Suppose that $\psi_1, \dots, \psi_m \in
\Sol(\Phi)$ are linearly independent over $F$. Then they are
linearly independent over $\bL$.
\end{lemma}

\begin{proof}
The proof is in the same spirit as the one for
Lemma~\ref{L:Bindependence}, and we omit it.
\end{proof}

\begin{corollary} \label{C:SolDimBound}
Let $\Phi \in \Mat_r(K)$.  Then $\Sol(\Phi)$ is an $F$-vector
space of dimension at most $r$.
\end{corollary}

\subsubsection{Fundamental matrix of solutions}
Given $\Phi \in \GL_r(K)$, suppose $\Psi \in \GL_r(L)$ satisfies
\[
  \sigma(\Psi) = \Phi\Psi.
\]
Then by Lemma~\ref{L:SolIndependence} and
Corollary~\ref{C:SolDimBound}, the columns of $\Psi$ form an
$F$-basis for $\Sol(\Phi)$.  The matrix $\Psi$ is called a
\emph{fundamental matrix for $\Phi$}.  It is useful to note that
$\Psi' \in \GL_r(L)$ is another fundamental matrix for $\Phi$ if and
only if $\Psi^{-1}\Psi'$ is fixed by $\sigma$.  That is, if and only
if $\Psi' = \Psi\delta$ for some $\delta \in \GL_r(F)$.

\subsection{The difference Galois group} \label{SS:GroupGamma}
Throughout this section we fix $\Phi \in \GL_r(K)$ and suppose that
$\Psi \in \GL_r(L)$ is a fundamental matrix for $\Phi$ with respect
to our $\sigma$-admissible triple $(F,K,L)$.

For a ring $R$, we let $\GL_{r/R}$ denote the $R$-group scheme of
$r\times r$ invertible matrices.  Its coordinate ring is $R[X,1/\det
X]$, where $X = (X_{ij})$ is an $r\times r$ matrix of independent
variables.  If $S$ is an $R$-algebra, we will as usual let $GL_r(S)$
denote the group of $S$-rational points on $\GL_{r/R}$.  For any
$R$-scheme $Z$, we let $Z_S := S \times_R Z$ be its base extension
to an $S$-scheme.

\subsubsection{Construction of $\Gamma$}
We define a $K$-algebra map $\nu : K[X,1/\det X] \to L$ by setting
$\nu(X_{ij}) := \Psi_{ij}$.  We let
\[
  \fp := \ker \nu, \qquad \Sigma := \im \nu
  = K[\Psi,1/\det \Psi] \subseteq L.
\]
We let $\Lambda$ be the fraction field of $\Sigma$.  Finally, we let
$Z = \Spec \Sigma$.  In this way $Z$ is the small closed subscheme
of $\GL_{r/K}$ such that $\Psi \in Z(L)$.

Now set $\Psi_1, \Psi_2 \in \GL_r(L\otimes_K L)$ to be the matrices
such that $(\Psi_1)_{ij} = \Psi_{ij} \otimes 1$ and $(\Psi_2)_{ij} =
1 \otimes \Psi_{ij}$, and let $\tPsi := \Psi_1^{-1}\Psi_2 \in
\GL_r(L\otimes_K L)$.  We define an $F$-algebra map $\mu :
F[X,1/\det X] \to L \otimes_K L$ by $\mu(X_{ij}) = \tPsi_{ij}$.  We
let
\[
  \fq := \ker \mu, \qquad \Delta := \im \mu,
\]
and finally we set $\Gamma = \Spec \Delta$.  In this way $\Gamma$ is
the smallest closed subscheme of $\GL_{r/F}$ such that $\tPsi \in
\Gamma(L \otimes_{K} L)$.  If we wish denote the dependence on
$\Psi$, we will write $Z_\Psi$ and $\Gamma_\Psi$ for these
spaces.

Among other things, we will see in Theorem~\ref{T:GammaGroup} that
$\Gamma$ is a closed subgroup of $\GL_{r/F}$ and that $Z$ is a
$\Gamma_K$-torsor under right-multiplication.

\subsubsection{The automorphisms $\bsigma_0$ and $\bsigma_1$}
\label{SSS:bsigma0} We define a natural $\sigma$-linear
automorphisms
\[
\bsigma_0, \bsigma_1 : L[X,1/\det X] \to L[X,1/\det X],
\]
by setting $\bsigma_0 X := X$ and $\bsigma_1 X := \Phi X$.  We note
that $\bsigma_0$ restricts to an automorphism of $R[X,1/\det X]$ for
any $F$-subalgebra $R$ of $L$, and that $\bsigma_1$ induces
automorphisms of $K[X,1/\det X]$ and $\Sigma[X,1/\det X]$.  We see
that
\[
  \bsigma_0 \fq = \fq, \qquad \bsigma_1 \fp = \fp.
\]
The first equality is clear since $\fq \subseteq F[X,1/\det X]$. For
the second, we note that for $h(X) \in K[X,1/\det X]$, we have
$\bsigma_1(h)(X) = \bsigma_0(h)(\Phi X)$, and so
\[
\bsigma_1(h)(\Psi) = \bsigma_0(h)(\Phi\Psi) =
\bsigma_0(h)(\sigma(\Psi)) = \sigma(h(\Psi)).
\]
Thus, $\bsigma_1 \fp = \fp$.  This equality implies further that
$\nu\bsigma_1 = \sigma\nu$.  The following lemma provides a
correspondence between the contraction and extension of ideals in
$L[X,1/\det X]$. See also \cite[Lem.~1.23]{vdpsing03}.

\begin{lemma} \label{L:SigmaIdealCorrespondence}
The functions between sets of ideals,
\[
\begin{array}{ccc}
\{ \fa \subseteq F[X,1/\det X] \} & \longleftrightarrow & \{ \fb
\subseteq L[X,1/\det X] \mid \bsigma_0\fb = \fb \}, \\
\fa & \to & (\fa) \\
\fb \cap F[X,1/\det X] & \gets & \fb
\end{array}
\]
are bijections.
\end{lemma}

\begin{proof}
Since $\bsigma_0$ is trivial on $F[X,1/\det X] \subseteq L[X,1/\det
X]$, these maps are well-defined.  One knows already that $(\fa)
\cap F[X,1/\det X] = \fa$ for all ideals $\fa \subseteq F[X,1/\det
X]$ (see \cite[\S VII.11]{zarsam2}).  Now let $\fb \subseteq
L[X,1/\det X]$ be an ideal with $\bsigma_0 \fb = \fb$, and let $\fa
:= \fb \cap F[X,1/\det X]$.  Letting $\{g_i\}_{i \in I}$ be an
$F$-basis of $F[X,1/\det X]$, we have that $\{g_i\}_{i\in I}$ is an
$L$-basis of $L[X,1/\det X]$.  For $h \in \fb$ we write $h = \sum
b_i g_i$, $b_i \in L$, and we let $l(h)$ be the number of $i \in I$
for which $b_i \neq 0$.  We show that $h \in (\fa)$ by induction on
$l(h)$.  If $l(h) = 0$ the result is clear.  If $l(h)=1$, then $h =
bg$ for some $b \in L^\times$ and $g \in \{ g_i \}$.  Moreover, then
$g \in \fa$.  Now suppose that $l(h) > 1$.  By multiplying by an
element of $L$ we can assume that $b_{i_1} = 1$ and that $b_{i_2}
\in L \setminus F$ for some $i_1, i_2 \in I$.  (If all $b_i \in F$,
then $h \in \fa$.)  One sees that
\[
  l(\bsigma_0 h - h) < l(h),
\]
and since $\bsigma_0 \fb = \fb$, we have $\bsigma_0 h - h \in \fb$.
Therefore, $\bsigma_0(h) - h \in (\fa)$.  Similarly,
$\bsigma_0(b_{i_2}^{-1}h) - b_{i_2}^{-1}h \in (\fa)$.  However,
\[
  (\sigma(b_{i_2}^{-1}) - b_{i_2}^{-1})h
  = \bigl(\bsigma_0(b_{i_2}^{-1}h) - b_{i_2}^{-1}h\bigr)-
  \sigma(b_{i_2}^{-1})\cdot(\bsigma_0 h -h).
\]
The left-hand side is non-zero, and the right-hand side is in
$(\fa)$. Therefore $h \in (\fa)$.
\end{proof}

\begin{proposition} \label{P:ZIsomGammaOverL}
Define a morphism of affine $L$-schemes $\phi := Z_L \to \GL_{r/L}$
so that on points $u \mapsto \Psi^{-1}u$ for $u \in
Z(\overline{L})$. Then $\phi$ factors through an isomorphism $\phi'
: Z_L \to \Gamma_L$ of affine $L$-schemes.
\end{proposition}

\begin{proof}
For commutative rings $R \subseteq S$ and for any ideal $I$ in
$R[X,1/\det X]$, we let $I_S$ denote its extension to $S[X,1/\det
X]$. Now the ideal $\fp \subseteq K[X,1/\det X]$ is the defining
ideal of the $K$-scheme $Z$, and $\fq \subseteq F[X,1/\det X]$ is
the defining ideal of of the $F$-scheme $\Gamma$. If we set
\[
  \alpha : L[X,1/\det X] \to L[X,1/\det X],
\]
to be the $L$-algebra homomorphism determined by setting $\alpha(X)
= \Psi^{-1}X$, then the map
\[
  \oalpha : L[X,1/\det X] \to L[X,1/\det X]\bigm/ \fp_L,
\]
induced by $\alpha$, is the map $\phi$ on the level of coordinate
rings.  It then suffices to prove that $\fq_L = \alpha^{-1}\fp_L$.

As noted in \S\ref{SSS:bsigma0}, we have that $\bsigma_0 \fq_L =
\fq_L$ and $\bsigma_1 \fp_L = \fp_L$. Furthermore,
\[
  \bsigma_1 \alpha X = \bsigma_1(\Psi^{-1} X) = (\sigma \Psi)^{-1}
  (\bsigma_1 X) = \Psi^{-1}X = \alpha\bsigma_0 X,
\]
and so $\bsigma_1\alpha = \alpha\bsigma_0$, which implies that
\[
  \bsigma_0 \alpha^{-1} \fp_L = \alpha^{-1}\fp_L.
\]
By Lemma~\ref{L:SigmaIdealCorrespondence}, it follows that
$\alpha^{-1} \fp_L$ is generated by $\alpha^{-1} \fp_L \cap
F[X,1/\det X]$.

Now we regard $L\otimes_K L$ as an $L$-algebra through the map $f
\mapsto f \otimes 1$.  If we let $\tmu : L[X,1/\det X] \to L
\otimes_K L$ be the unique $L$-algebra homomorphism such that $\tmu
X = \Psi_2$, then we note that the composition
\[
  F[X,1/\det X] \stackrel{\alpha}{\to} L[X,1/\det X]
  \stackrel{\tmu}{\to} L \otimes_K L
\]
is in fact $\mu$.  Since $L$ is a field, the map $L[X,1/\det
X]/\fp_L \to L\otimes_K L$ induced by $\tmu$ is injective.
Therefore,
\[
  \fq = \alpha^{-1}\fp_L \cap F[X,1/\det X],
\]
and by our argument in the previous paragraph, $\fq_L =
\alpha^{-1}\fp_L$.
\end{proof}

\begin{corollary} \label{C:SigmaSimple}
The ideal $\fp \subseteq K[X,1/\det X]$ is maximal among proper
$\bsigma_1$-invariant ideals.
\end{corollary}

\begin{proof}
Let $\fm \supseteq \fp$ be a proper ideal of $K[X,1/\det X]$ such
that $\bsigma_1 \fm \subseteq \fm$.  Because $K[X,1/\det X]$ is
noetherian, it follows that $\bsigma_1 \fm = \fm$.  Now $\alpha^{-1}
\fm_L \supseteq \alpha^{-1} \fp_L = \fq_L$, and we see easily that
$\bsigma_0 \alpha^{-1}\fm_L = \alpha^{-1}\fm_L$. Therefore, by
Lemma~\ref{L:SigmaIdealCorrespondence},
\[
  \alpha^{-1}\fm_L = (\alpha^{-1}\fm_L \cap F[X,1/\det X])_L.
\]
Let $\fa \subseteq F[X,1/\det X]$ be a maximal ideal that contains
$\alpha^{-1}\fm_L \cap F[X,1/\det X]$, and let $E := F[X,1/\det
X]/\fa$, which is a finite extension of $F$. By
Lemma~\ref{L:SigmaIdealCorrespondence}, we see that $\fa = \fa_L
\cap F[X,1/\det X]$, and it follows that there is an isomorphism
$\beta : L[X,1/\det X]/\fa_L \iso L \otimes_F E$.  Now if we
consider the maps
\[
  \Pi : L[X,1/\det X] \stackrel{\alpha^{-1}}{\to} L[X,1/\det X]
  \stackrel{\beta}{\to} L \otimes_F E,
\]
we see that $\fm_L \subseteq \ker \Pi$. If we let $\pi : K[X,1/\det
X] \to L \otimes_F E$ be the restriction of $\Pi$, then easily $\fm
\subseteq \ker \pi$ and $\ker \pi$ is a proper ideal.  Moreover,
since $\alpha^{-1}\bsigma_1 = \bsigma_0 \alpha^{-1}$, it follows
that $\ker \pi$ is a $\bsigma_1$-invariant ideal of $K[X,1/\det X]$.
Therefore, the maximality of $\fm$ implies that $\fm = \ker \pi$.

Now let $\Psi' \in \GL_r(L\otimes_F E)$ be defined by $\Psi'_{ij} =
\pi(X_{ij})$.  The automorphism $\sigma$ on $L$ extends to an
automorphism of $L\otimes_F E$  by acting by the identity on $E$,
and it is easily seen that $(L \otimes_F E)^\sigma = E$. In this way
$\sigma(\Psi') = \Phi\Psi'$, and this implies that the matrix
$\delta := (\Psi')^{-1}\Psi \in \GL_r(E)$.  Now $\delta$ induces an
automorphism on $(K \otimes_F E)[X,1/\det X]$ via
\[
  \delta \cdot h(X) := h(X\delta).
\]
If we extend $\pi$ to $\pi' : (K \otimes_F E)[X,1/\det X] \to L
\otimes_F E$ by the identity on $E$, then we see that we have the
extended ideals
\[
  \fp_{K \otimes_F E} = \fp \otimes_F E \subseteq
  \fm \otimes_F E \subseteq \ker \pi' = \delta \cdot (\fp \otimes_F E).
\]
But $(K \otimes_F E)[X,1/\det X]$ is a noetherian ring, and so $\fp
\otimes_F E \subseteq \delta \cdot (\fp \otimes_F E)$ implies that
$\fp \otimes_F E = \delta \cdot (\fp \otimes_F E)$.  Thus,
\[
\fp \otimes_F E = \fm \otimes_F E = \delta \cdot (\fp \otimes_F E).
\]
Now $(K \otimes_F E)[X,1/\det X]$ is a free $K[X,1/\det X]$-module,
since $E$ is a vector space over $F$, and is therefore faithfully
flat over $K[X,1/\det X]$.  So by intersecting with $K[X,1/\det X]$
we see that $\fp = \fm$ (see \cite[Thm.~7.5]{mats2}).
\end{proof}

\subsubsection{Contracted ideals of $\Sigma[X,1/\det X]$}
The following lemma is a companion to
Lemma~\ref{L:SigmaIdealCorrespondence}, and relies on the preceding
corollary.  See also \cite[Lem.~1.11]{vdpsing97}.

\begin{lemma} \label{L:SigmaXIdeals}
Let $\fb \subset \Sigma[X,1/\det X]$ be an ideal that is
$\bsigma_0$-invariant. Then $\fb$ is generated by $\fb \cap
F[X,1/\det X]$.
\end{lemma}

\begin{proof}
Let $\fa := \fb \cap F[X,1/\det X]$.  Let $\{ g_i \}_{i\in I}$ be an
$F$-basis for $F[X,1/\det X]$ such that $I = I_\fa \cup I_1$, where
$\{g_i \}_{i\in I_\fa}$ is an $F$-basis for $\fa$.  Choose a subset
$J \subseteq I_1$ minimal so that $\fb \cap \sum_{i \in J} \Sigma
g_i$ contains a non-zero element of $\fb/(\fa)$. Pick $j \in J$, and
let
\[
  \fm := \{ b \in \Sigma \mid {\textstyle \exists \sum_{i \in J} b_i g_i \in
  \fb/(\fa), b_j = b}\}
\]
Since $\bsigma_0 \fb = \fb$, and since each $g_i$ is fixed by
$\bsigma_0$, it follows that $\fm$ is a non-zero $\sigma$-invariant
ideal of $\Sigma$.  However, by Corollary~\ref{C:SigmaSimple},
$\Sigma$ has no $\sigma$-invariant ideals other than $\{ 0 \}$ and
$\Sigma$. Thus, $\fm = \Sigma$.  Therefore there exists $h \in \fb$
so that $h = \sum_{i \in J} b_ig_i \mod (\fa)$ and $b_j = 1$.  Now
$\bsigma_0(h) - h$ is supported on a proper subset of $J$ modulo
$(\fa)$, and so it must be $0$ modulo $(\fa)$ by the minimality of
$J$. Therefore each $b_i \in F$, and thus $\sum_{i \in J} b_ig_i \in
(\fa)$, which is a contradiction.
\end{proof}

\begin{proposition} \label{P:ZGammaKProduct}
Define a morphism of affine $K$-schemes $\psi : Z \times Z \to Z
\times \GL_{r/K}$ so that on points $(u,v) \mapsto (u,u^{-1}v)$ for
$u$, $v \in Z(\oK)$.  Then $\psi$ factors through an isomorphism $Z
\times Z \to Z \times \Gamma_K$ of affine $K$-schemes.
\end{proposition}

\begin{proof}
Again we work on the level of coordinate rings and maintain
conventions and definitions in the proof of
Proposition~\ref{P:ZIsomGammaOverL}.  The ring $\Sigma \subseteq L$
is isomorphic to the coordinate ring of $Z$ over $K$.  Likewise, the
ring $\Sigma[X,1/\det X]/\fp_\Sigma$ is the coordinate ring of $Z
\times Z$, and $\Sigma[X,1/\det X]/\fq_\Sigma$ is the coordinate
ring of $Z \times \Gamma_K$.  The $L$-algebra automorphism $\alpha$
in the proof of the previous proposition restricts to an
automorphism of $\Sigma[X,1/\det X]$, and in this way the
homomorphism
\[
  \oalpha : \Sigma[X,1/\det X] \to \Sigma[X, 1/\det X] \bigm/
  \fp_\Sigma
\]
induced by $\alpha$ represents the morphism $\psi$ of affine
$K$-schemes.  We then need to show that $\fq_\Sigma =
\alpha^{-1}\fp_\Sigma$.

Let $\fa = \alpha^{-1}\fp_\Sigma \cap F[X,1/\det X]$.  By
Lemma~\ref{L:SigmaXIdeals},
\[
\alpha^{-1}\fp_{\Sigma} = \fa_{\Sigma}.
\]
Then as in the proof of Proposition~\ref{P:ZIsomGammaOverL}, $\fq_L
= \alpha^{-1}\fp_L$ and both are now equal to $\fa_L$. Since $\fq$
and $\fa$ are both ideals in $F[X,1/\det X]$ and $F$ and $L$  are
fields, it follows that $\fq = \fa$.
\end{proof}

\begin{lemma} \label{L:SubgroupLemma}
Let $G$ be a group, and let $A$ and $B$ be subsets of $G$, $A$
non-empty, such that the map
\[
  (u,v) \mapsto (u,u^{-1}v) : A \times A \to A \times G
\]
factors through a bijection $\phi : A \times A \to A \times B$.
Then $B$ is a subgroup of $G$ and $A$ is stable under
right-multiplication by elements of $B$.  Moreover, under the action
of $B$ by right-multiplication, $A$ becomes a principal homogeneous
space for $B$.
\end{lemma}

\begin{proof}
This is a simple exercise.
\end{proof}

\subsubsection{The Galois group $\Gamma$}
The previous propositions and lemmas culminate in the following
theorem, saying that $\Gamma$ is in fact an affine group scheme.  We
call the group $\Gamma$, or $\Gamma_\Psi$ if we wish to recall the
dependence on $\Psi$, the \emph{Galois group of the system
$\sigma(\Psi) = \Phi\Psi$}.

\begin{theorem} \label{T:GammaGroup}
Let $(F,K,L)$ be a $\sigma$-admissible triple for an automorphism
$\sigma : L \to L$.  Suppose we have $\Phi \in GL_r(K)$ and $\Psi
\in GL_r(L)$ so that $\sigma(\Psi) = \Phi\Psi$.  Then $\Gamma :=
\Gamma_\Psi$ is a closed $F$-subgroup scheme of $\GL_{r/F}$, and the
closed $K$-subscheme $Z := Z_\Psi$ of $\GL_{r/K}$ is stable under
right-multiplication by $\Gamma_K$ and is a $\Gamma_K$-torsor.
\end{theorem}

\begin{proof}
Since $Z(\Sigma)$ is non-empty, Propositions~\ref{P:ZIsomGammaOverL}
and~\ref{P:ZGammaKProduct} imply that $(u,v) \mapsto (u,u^{-1}v):
Z(\Sigma) \times Z(\Sigma) \to Z(\Sigma) \times \Gamma(\Sigma)$ is a
bijection.  Lemma~\ref{L:SubgroupLemma} and the Yoneda lemma
\cite[\S 1.2--1.4]{water} imply that $\Gamma_\Sigma$ is a subgroup
of $\GL_{r/\Sigma}$ and that $Z_\Sigma$ is a $\Gamma_\Sigma$-torsor.
Since the inclusion $F \to \Sigma$ is faithfully flat, we see that
$\Gamma$ is a closed $F$-subgroup scheme of $\GL_{r/F}$ by flat
descent \cite[\S 17.1--17.3]{water}. Similarly, since the inclusion
$K \to \Sigma$ is faithfully flat, $Z$ admits the structure of a
$\Gamma_K$-torsor.
\end{proof}

\subsection{Criterion for smoothness}

We continue with the notation of the previous section, and in
particular have fixed a $\sigma$-admissible triple $(F,K,L)$
together with $\Phi \in \GL_r(K)$, $\Psi \in \GL_r(L)$ satisfying
$\sigma(\Psi) = \Phi\Psi$.  In this section, we explore when
$\Gamma$ is smooth over $\oF$, that is, when the coordinate ring of
$\Gamma_{\oF}$ is reduced.

\begin{theorem} \label{T:Smoothness}
Suppose $K$ is algebraically closed in the fraction field $\Lambda$
of $\Sigma$. Then
\begin{enumerate}
\item[(a)] The $K$-scheme $Z$ is absolutely irreducible and is smooth over
$\oK$.
\item[(b)] The $F$-scheme $\Gamma$ is absolutely irreducible and is smooth over
$\oF$.
\item[(c)] The dimension of $\Gamma$ over $F$ is equal to the
transcendence degree of $\Lambda$ over $K$.
\end{enumerate}
\end{theorem}

\begin{proof}
The ideal $\fp \subseteq K[X,1/\det X]$ is prime.  The field
$\Lambda$ is separable over $K$, since it is a subfield of $L$. That
$K$ is algebraically closed in $\Lambda$ then implies that
$\fp_{\oK} \subseteq \oK[X,1/\det X]$ is prime \cite[VII.11,
Thm.~39]{zarsam2}. Thus $Z$ is smooth over $\oK$. Because
$\Gamma_{\oK} \cong Z_{\oK}$, $\Gamma$ must be smooth over $\oF$.
By construction the transcendence degree of $\Lambda$ over $K$ is
equal to the dimension of $Z$, which is equal to the dimension
of~$\Gamma$.
\end{proof}

\subsubsection{The case $(\FF_q(t),\ok(t),\LL)$}
This case is of particular interest to our applications to
$t$-motives.  It turns out that in this case, all Galois groups are
smooth, via the following proposition.  We continue with our usual
notation.

\begin{proposition}
Suppose $(F,K,L) = (\FF_q(t),\ok(t),\LL)$, and suppose that $\Phi
\in \GL_r(\ok(t))$ and $\Psi \in \GL_r(\LL)$ satisfy $\Psi^{(-1)} =
\Phi\Psi$.  Then $\ok(t)$ is algebraically closed in $\Lambda =
\ok(t)(\Psi)$.
\end{proposition}

\begin{proof}
Let $f \in \Lambda \cap \okt$, and consider the field $H :=
\ok\bigl(t;\ f^{(i)} : i \in \ZZ \bigr)$ obtained by adjoining all
of the twists of $f$ to $\ok(t)$.  Each $f^{(i)}$ is algebraic over
$\ok(t)$, and so $H/\ok(t)$ is algebraic. Since $\Lambda$ is
finitely generated as a field over $\ok(t)$, so is $H$.  Thus
$[H:\ok(t)] < \infty$. Furthermore, $H$ is invariant under $\sigma$
and $\sigma^{-1}$.

The field $H$ is the function field of a smooth projective curve $X$
over $\ok$, and the inclusion $\ok(t) \subseteq H$ provides a
surjective morphism $X \to \PP^1_{\ok}$ over $\ok$.  Now $\sigma: H
\to H$ induces an automorphism $\tau : X \to X$ as a scheme over
$\FF_q$.  Because $\sigma$ leaves the integral closure of $\ok[t]$
in $H$ invariant, the points $\infty_1, \dots, \infty_d$ in $X$
above the point $\infty$ in $\PP^1_{\ok}$ are permuted by $\sigma$.
Thus we can construct an effective divisor $I$ of $X$ such that
$\tau(I)=I$ and $\Supp(I)=\{\infty_1, \dots, \infty_d\}$.  Now for
$N \geq 1$ sufficiently large, the field $H$ is generated over
$\ok(t)$ by the functions in the finite dimensional $\ok$-vector
space
\[
  S := \Gamma(X,N\cdot I) \subseteq H.
\]
By our assumptions on $I$, this space is invariant under $\sigma$
and $\sigma^{-1}$. If the entries of $\bff := [f_1, \dots,
f_m]^{\tr}$ form a $\ok$-basis for $S$, then there is a matrix $A
\in \GL_m(\ok)$ so that $\sigma(\bff) = A\bff$.  If $\bg \in
\Mat_{m\times 1}(S)$ and $\bg = B\bff$ for some $B \in \GL_m(\ok)$,
then
\[
\sigma(\bg) = B^{(-1)} A B^{-1}  \bg.
\]
By the theory of Lang isogenies \cite{lang56}, we can pick a $B \in
\GL_m(\ok)$ so that
\[
  B^{-1}B^{(1)} = A^{(1)},
\]
and if we let $\bg := B\bff$, then $\sigma(\bg) = \bg$. Thus $S$
contains a $\ok$-basis $\bg$ that is fixed by $\sigma$, and $H =
\ok(t,\bg)$.  Let $g$ be an entry of $\bg$. Then $g \in \okt \cap
\LL^\sigma = \FF_q(t)$. Thus $[H:\ok(t)] = 1$.
\end{proof}

\begin{corollary}
Let $(F,K,L) = (\FF_q(t),\ok(t),\LL)$, and suppose that $\Phi \in
\GL_r(\ok(t))$ and $\Psi \in \GL_r(\LL)$ satisfy $\sigma(\Psi) =
\Phi\Psi$.  Then the Galois group $\Gamma$ of $\Psi$ is smooth over
$\oFqt$.
\end{corollary}

\subsection{The Galois action} \label{SS:GaloisAction}

In this section we will assume that $K$ is algebraically closed in
$\Lambda$, and so in particular by Theorem~\ref{T:Smoothness},
$\Gamma_{\oF}$ and $Z_{\oK}$ are reduced and irreducible.

\subsubsection{$\sigma$-automorphisms of $\Sigma$ and $\Lambda$}
\label{SSS:sigmaautos} Let $\Aut_{\sigma}(\Sigma/K)$ denote the
group of automorphisms of $\Sigma$ over $K$ that commute with
$\sigma$. Similarly we define $\Aut_\sigma(\Lambda/K)$.  In fact, it
is true that
\[
  \Aut_\sigma(\Sigma/K) = \Aut_\sigma(\Lambda/K).
\]
Indeed every automorphism $\xi \in \Aut_\sigma(\Sigma/K)$ extends
uniquely to an automorphism in $\Aut_\sigma(\Lambda/K)$.  On the
other hand, if $\eta \in \Aut_\sigma(\Lambda/K)$, then as matrices
in $\Mat_r(L)$, $\sigma(\eta(\Psi)) = \eta(\sigma(\Psi)) =
\eta(\Phi\Psi) = \Phi\eta(\Psi)$. Thus, $\eta(\Psi) = \Psi\gamma$
for some $\gamma \in \GL_r(F)$, and so $\eta(\Psi) \in
\Mat_r(\Sigma)$.  Therefore $\eta$ restricted to $\Sigma$ takes
values in $\Sigma$.

\subsubsection{The action of $\Gamma(F)$}
For $\gamma \in \Gamma(F)$, we have an automorphism of $K$-schemes
$\gamma : Z \to Z$ defined by right multiplication by $\gamma$.  On
the level of coordinate rings, the induced map is
\[
  \gamma = (h(X) \mapsto h(X\gamma)) : \Sigma \to \Sigma,
\]
which is a $K$-linear automorphism that commutes with the action of
$\sigma$.  Thus we have a group homomorphism,
\[
  \kappa: \Gamma(F) \to \Aut_\sigma(\Lambda/K),
\]
which is easily seen to be injective.  Now if $\delta \in
\Aut_\sigma(\Lambda/K)$, then $\delta$ induces an automorphism of
the non-empty $Z(\Lambda)$ that is right-multiplication by an
element of $\gamma \in \Gamma(\Lambda)$. That $\delta$ commutes with
$\sigma$ implies that $\gamma \in \Gamma(\Lambda^\sigma) =
\Gamma(F)$.  Thus $\kappa$ is an isomorphism.

\subsubsection{Base extensions} \label{SSS:sigmaBaseExtensions}
Given our $\sigma$-admissible triple
$(F,K,L)$, we choose an extension of $\sigma$ to an automorphism of
$\oL$. Then $\oL^\sigma$ is an algebraic extension of $F$.  Indeed,
the monic irreducible polynomial of any $h \in \oL^\sigma$ over $L$
must have coefficients in $L^\sigma = F$.  Thus if we let $\bF =
\oL^\sigma$, then $(\bF, \oK, \oL)$ is a $\sigma$-admissible triple.
The Galois group $\Gamma'$ defined by the system $\sigma(\Psi) =
\Phi\Psi$ defined with respect to $(\bF, \oK, \oL)$ is seen to be
$\Gamma_{\bF}$ by Propositions~\ref{P:ZIsomGammaOverL}
and~\ref{P:ZGammaKProduct}. If we let $\tSigma$ be the coordinate of
$Z_{\oK}$ and $\tLambda$ be its fraction field, then we see that
\[
  \Gamma(\bF) \cong \Aut_{\sigma}(\tSigma/\oK) =
  \Aut_{\sigma}(\tLambda/\oK).
\]

Furthermore, for $n\geq 1$, let $\bF_n = \oL^{\sigma^n}$, and
suppose that $(\bF_n,\oK,\oL)$ is $\sigma^n$-admissible.  Then
$\Psi$ is a fundamental matrix for $\Phi_n :=
\sigma^{n-1}(\Phi)\cdots\sigma(\Phi)\Phi$.  Again by
Propositions~\ref{P:ZIsomGammaOverL} and~\ref{P:ZGammaKProduct}, we
see that the Galois group of this system of equations is
$\Gamma_{\bF_n}$.  And thus,
\[
  \Gamma(\bF_n) \cong \Aut_{\sigma^n}(\tSigma/\oK) =
  \Aut_{\sigma^n}(\tLambda/\oK).
\]

\subsubsection{Galois action for $\Gamma(\oF)$}
Continuing with the notation of the previous paragraphs, suppose
that $\oF = \cup \bF_n$.  Then every element of $\Gamma(\oF)$
induces an automorphism of $\tLambda/\oK$ that commutes with
$\sigma^n$ for all $n \gg 0$.  In this case, we will call this the
induced action of $\Gamma(\oF)$ on $\tLambda$.

\subsubsection{The case $(\FF_q(t),\ok(t),\LL)$}  It is worth pointing out
that the situation is quite nice in our usual setting.  For $n\geq
1$, the triple $(\FF_{q^n}(t),\ok(t),\LL)$ is $\sigma^n$-admissible.
As in \S\ref{SSS:Twisting}, there is a canonical extension of
$\sigma$ to $\genpow{\KK}{t} \supseteq \oLL$. Furthermore, we see
that
\[
  \oLL^{\sigma^n} = \okt^{\sigma^n} = \oFqt^{\sigma^n} =: \bF_n,
\]
and so $(\bF_n, \okt, \oLL)$ is a $\sigma^n$-admissible triple.
Every element of $\oFqt$ is fixed by some power of $\sigma$, and so
\[
  \oFqt = \bigcup_{n\geq 1} \bF_n.
\]
We now return to the general situation, but it is important to note
that the following theorem applies to Galois groups in the usual
$(\FF_q(t),\ok(t),\LL)$ setting.

\begin{theorem} \label{T:GammaFixedIsoK}
Let $\Phi \in \GL_r(K)$, and suppose that $\Psi \in \GL_r(L)$ is a
fundamental matrix for $\Phi$.  Assume that $K$ is algebraically
closed in $\Lambda = K(\Psi)$.  Fix an extension of $\sigma$ to
$\oL$, and let $\bF_n := \oL^{\sigma^n}$. Suppose that
$(\bF_n,\oK,\oL)$ is $\sigma^n$-admissible for each $n \geq 1$, and
suppose that $\oF = \cup \bF_n$.  Let $\tSigma$ be the coordinate
ring of $Z_{\oK}$ and let $\tLambda$ be its fraction field, both
considered subrings of $\oL$.
\begin{enumerate}
\item[(a)] The subfield of $\tLambda$ fixed by $\Gamma(\oF)$
is $\oK$.
\item[(b)] The elements of $\Lambda$ fixed by
$\Gamma(\oF)$ are precisely $K$.
\end{enumerate}
\end{theorem}

\begin{proof}
See \cite[Lem.~1.28]{vdpsing97}.  Suppose $f \in \tLambda$ is fixed
by $\Gamma(\oF)$.  We consider $f \in \tLambda$ to be a function $f
: Z_{\oK} \to \PP^1_{\oK}$.  For $i=1, 2$, we consider the two maps
of $\oK$-schemes
\[
  g_i : Z_{\oK} \times \Gamma_{\oK} \to Z_{\oK} \times
  Z_{\oK} \stackrel{\pi_i}{\to} Z_{\oK} \stackrel{f}{\to}
  \PP^1_{\oK},
\]
where $\pi_i$ is the $i$-th projection.  Because $f$ is
$\Gamma(\oF)$-invariant and because $\Gamma(\oF)$ is dense in
$\Gamma_{\oK}$ since $\Gamma$ is smooth over $\oF$, we must have
$g_1=g_2$. Therefore, $f \circ \pi_1 = f \circ \pi_2$, which implies
that $f$ is constant. This proves part (a). Part (b) follows from
part (a) and the assumption that $K$ is algebraically closed in
$\Lambda$.
\end{proof}

\subsubsection{Remark}
If $\Gamma(F)$ is Zariski dense in $\Gamma$, then it follows that
$\Lambda^{\Gamma(F)} = \ok(t)$.

\subsection{The group $\Gamma$ and $t$-motives}
\label{SS:GammaPsiAndtMotives} Given a $t$-motive $M$, we defined
the Galois group $\Gamma_M$ of $M$ in
\S\ref{SSS:tMotiveGaloisGroups}.  Associated to $M$ we can also
choose a matrix $\Phi \in \GL_r(\ok(t))$ that represents
multiplication by $\bsigma$ on $M$.  Let $\Psi \in \GL_r(\LL)$ be a
rigid analytic trivialization of $\Phi$.  We will show that
$\Gamma_M$ is isomorphic to $\Gamma := \Gamma_\Psi$ over $\FF_q(t)$.

\subsubsection{$t$-motives and $\sigma$-semilinear equations}
Let $M$ be a $t$-motive.  We fix the following notation throughout
this section.  Let $\bm$ be a basis for $M$, and let $\Phi \in
\GL_r(\ok(t))$ represent multiplication by $\bsigma$ on $M$. We pick
a rigid analytic trivialization $\Psi \in \GL_r(\LL)$ for $M$, which
is at the same time a fundamental matrix for $\Phi$.

Let $M^u_v := M^{\otimes u} \otimes (M^{\vee})^{\otimes v}$. Because
$\cT_M$ is Tannakian, if $N$ is any $t$-motive in $\cT_M$, then $N$
is the subquotient of a direct sum of various $M^u_v$, and vice
versa.  It follows from Propositions~\ref{P:RATrivialization}(b)
and~\ref{P:RATKerQuot}(b) that the entries of a fundamental matrix
$\Psi_N$ for $N$ are in $\Sigma$, and in fact we can take $\Psi_N
\in \GL_s(\Sigma)$ for some $s$.

For an $\FF_q(t)$-algebra $R$, we let $\Sigma^{(R)} := R
\otimes_{\FF_q(t)} \Sigma$.

\begin{lemma} \label{L:RATBijectiveForSigma}
For any $t$-motive $N$ in $\cT_M$ and $\FF_q(t)$-algebra $R$, the
natural map,
\[
  \Sigma^{(R)} \otimes_{\FF_q(t)} N^{\rB} \to
  \Sigma^{(R)} \otimes_{\ok(t)} N
\]
is bijective.
\end{lemma}

\begin{proof}
Let $\kappa$ be the map defined in the statement of the lemma. Thus
as above we can pick a basis $\bn$ for $N$ and a rigid analytic
trivialization $\Psi_N \in \GL_s(\Sigma)$ with respect to $\bn$. By
Proposition~\ref{P:RATrivialization}(b), $\Psi_N^{-1} \bn$ is an
$\FF_q(t)$-basis for $N^{\rB}$.  Now $1 \otimes (\Psi_N^{-1}\bn)$ is
a $\Sigma^{(R)}$-basis of $\Sigma^{(R)} \otimes_{\FF_q(t)} N^{\rB}$
(here and elsewhere $1 \otimes A$ for a matrix $A$ is the matrix of
the same dimension whose entries are each tensored by 1 on the
left). If $\bff \in \Mat_{1 \times s}(\Sigma^{(R)})$, then
\[
  \kappa \bigl( (\bff \otimes 1)\cdot(1 \otimes
  (\Psi_N^{-1}\bn)) \bigr) = (\bff\Psi_N^{-1} \otimes 1)\cdot(1
  \otimes \bn).
\]
The entries of $(\Psi_N^{-1} \otimes 1)\cdot(1 \otimes \bn)$ are
in the image of $\kappa$, and
\[
  \Psi_N \cdot (\Psi_N^{-1} \otimes 1) \cdot (1
  \otimes \bn) = 1 \otimes \bn.
\]
Thus $\kappa$ is surjective.  Since $\Psi_N \in \GL_s(\Sigma)$, the
map $\kappa$ is bijective.
\end{proof}

\begin{theorem} \label{T:tMotivestoReps}
Let $M$ be a $t$-motive, and let $N$ be a $t$-motive in $\cT_M$.
If we consider $N^{\rB}$ to be an algebraic group over $\FF_q(t)$,
then there is a natural representation
\[
  \xi_N : \Gamma \to \GL(N^{\rB})
\]
over $\FF_q(t)$ that is functorial in $N$.
\end{theorem}

\begin{proof}
Since every $t$-motive $N$ in $\cT_M$ is constructed from $M$ via
tensor products, duals, and subquotients, to define this
representation it suffices to define it on $M^{\rB}$ itself.
Functoriality in $N$ will be automatic.

To define the representation on $M^{\rB}$, it suffices by the Yoneda
lemma \cite[\S 1.2--1.4]{water} to define a representation
\[
  \xi_M^{(R)} : \Gamma(R) \to \GL(R \otimes_{\FF_q(t)} M^{\rB})
\]
for every $\FF_q(t)$-algebra $R$ and show that it is functorial in
$R$.  Let $R$ be an $\FF_q(t)$-algebra, and let $\gamma \in
\Gamma(R)$. Define
\[
  \Xi^{(R)}(\gamma) := \gamma \otimes 1 = (h(\Psi) \otimes m \mapsto
  h(\Psi\gamma) \otimes m) : \Sigma^{(R)} \otimes_{\ok(t)}
   M \to \Sigma^{(R)} \otimes_{\ok(t)} M,
\]
which is an isomorphism of $\ok(t)$-vector spaces.  Now by
Lemma~\ref{L:RATBijectiveForSigma}, $R\otimes_{\FF_q(t)} M^{\rB}$
spans $\Sigma^{(R)} \otimes_{\ok(t)} M$ as a $\Sigma^{(R)}$-module.
Let $\xi^{(R)}(\gamma)$ be the restriction of $\Xi^{(R)}(\gamma)$ to
$R\otimes_{\FF_q(t)} M^{\rB}$.

We claim that the image of $\xi^{(R)}(\gamma)$ is $R
\otimes_{\FF_q(t)} M^{\rB}$.  Indeed, since $\Psi^{-1}\bm$ forms an
$\FF_q(t)$-basis of $M^{\rB}$ by
Proposition~\ref{P:RATrivialization}(b), for $\bff \in \Mat_{1
\times r}(R)$, we have
\[
  \xi^{(R)}(\gamma) : \bff \cdot (1 \otimes \Psi^{-1})\bm \mapsto
  \bff \cdot \gamma^{-1}(1\otimes \Psi^{-1})\bm.
\]
Thus $\xi^{(R)}(\gamma)$ is an $R$-linear automorphism of $R
\otimes_{\FF_q(t)} M^{\rB}$.  It is straightforward to check that
this is construction is functorial in $R$, and so we have defined a
homomorphism $\xi_M : \Gamma \to \GL(M^{\rB})$.
\end{proof}

\begin{corollary} \label{C:MBFaithful}
Let $M$ be a $t$-motive.  The representation $\xi_M : \Gamma \to
\GL(M^{\rB})$ is faithful.
\end{corollary}

\begin{proof}
As defined in the the proof of the previous theorem we see easily
that $\xi^{(R)} : \Gamma(R) \to \GL(R \otimes_{\FF_q(t)} M^{\rB})$
is injective for all $\FF_q(t)$-algebras $R$.
\end{proof}

\subsubsection{The functor $\xi_M$}
For a $t$-motive $M$, if $\Phi \in \GL_r(\ok(t))$ represents
multiplication by $\bsigma$ on $M$ and if $\Psi \in \GL_r(\LL)$ is a
rigid analytic trivialization of $\Phi$, then
Theorem~\ref{T:tMotivestoReps} defines a functor
\[
  \xi_M : \cT_M \to \Rep{\Gamma}{\FF_q(t)}.
\]
It is straightforward to check that $\xi_M$ is a tensor functor.
Let
\[
\eta_M : \Rep{\Gamma_M}{\FF_q(t)} \iso \cT_M
\]
be the equivalence of categories defined in
\S\ref{SSS:tMotiveGaloisGroups}. Letting
\[
F : \Rep{\Gamma}{\FF_q(t)} \to \Vector{\FF_q(t)}
\]
be the forgetful functor, we see immediately that $\omega_M = F
\circ \xi_M$.  Thus by \cite[Cor.~II.2.9]{delmil82}, there is a
unique homomorphism $\pi_M : \Gamma \to \Gamma_M$ over $\FF_q(t)$ so
that the natural functor $\tau_M : \Rep{\Gamma_M}{\FF_q(t)} \to
\Rep{\Gamma}{\FF_q(t)}$ induced by $\pi_M$ satisfies
\[
  \xi_M \circ \eta_M = \tau_M.
\]

\begin{proposition} \label{P:xiMFullyFaithful}
Let $M$ be a $t$-motive.  Suppose that $\Phi \in \GL_r(\ok(t))$
represents multiplication by $\bsigma$ on $M$ and that $\Psi \in
\GL_r(\LL)$ is a rigid analytic trivialization for $\Phi$.  Then the
functor
\[
  \xi_M : \cT_M \to \Rep{\Gamma}{\FF_q(t)}
\]
is fully faithful.
\end{proposition}

\begin{proof}
For any $t$-motives $N$ and $P$ in $\cT_M$, there is a natural
isomorphism of $\FF_q(t)$-vector spaces, $\Hom_{\cT_M}(P,N) \cong
\Hom_{\cT_M}(\one,\Hom(P,N))$. Thus it suffices to prove full
faithfulness when $P = \one$.  Now $\Hom_{\cT_M}(\one,N) \cong N
\cap N^{\rB} = \{ n \in N \mid \bsigma n = n \}$, and this provides
an injection $\Hom_{\cT_M}(\one,N) \hookrightarrow \Hom_{\Gamma}
(\one^{\rB},N^{\rB})$.

Conversely suppose that $\phi : \one^{\rB} \to N^{\rB}$ is a
$\Gamma$-morphism.  Pick a $\ok(t)$-basis $\bn$ for $N$. Then
$\phi(1) = \bh(\Psi) \cdot \bn$ for some $\bh(\Psi) \in
\Mat_{1\times s}(\Sigma)$ by Lemma~\ref{L:RATBijectiveForSigma}. Let
$E/\FF_q(t)$ be a finite extension of fields.  We see that for
$\gamma \in \Gamma(E)$, the action of $\xi^{(E)}(\gamma) :=
\xi_M^{(E)}(N)(\gamma)$ on $E \otimes_{\FF_q(t)} N^{\rB}$ is simply
the restriction of the natural map $\Xi^{(E)}(\gamma) = 1 \otimes
\gamma: \Sigma^{(E)} \otimes_{\ok(t)} N \to \Sigma^{(E)}
\otimes_{\ok(t)} N$ to $E \otimes_{\FF_q(t)} N^{\rB}$. Since $\phi$
is a $\Gamma$-morphism, it follows that $\xi^{(E)}(\gamma)(\phi(1))
= \phi(1)$ for all $\gamma \in \Gamma(E)$.  Thus,
\[
  \bh(\Psi)\cdot \bn = \phi(1) = \xi^{(E)}(\gamma)(\phi(1)) =
  \bh(\Psi\gamma)\cdot \bn, \quad \gamma \in \Gamma(E).
\]
Because $\bn$ is a $\Sigma^{(E)}$-basis of $\Sigma^{(E)}
\otimes_{\ok(t)} N$, the entries of $\bh(\Psi)$ must each be fixed
by every $\gamma \in \Gamma(E)$.  By varying over all $E/\FF_q(t)$
finite, Theorem~\ref{T:GammaFixedIsoK} implies that $\bh(\Psi) \in
\Mat_{1 \times s}(\ok(t))$.  Thus $\phi(1) \in N \cap N^{\rB}$.
\end{proof}

\begin{lemma} \label{L:InvariantSubspace}
Let $\Phi \in \GL_r(\ok(t))$, and suppose that $\Psi \in \GL_r(\LL)$
is a fundamental matrix for $\Phi$.  Suppose that $W \subseteq
\Lambda^{\oplus s}$ is a vector subspace over $\Lambda$ such that
for every finite extension of fields $E/\FF_q(t)$,
\[
  \Gamma(E) \cdot (E \otimes_{\FF_q(t)} W) \subseteq E
  \otimes_{\FF_q(t)} W.
\]
Then $W$ has a system of defining equations over $\ok(t)$.
\end{lemma}

\begin{proof}
Suppose that $W$ has dimension $s-m$, and let $A(\Psi) \in
\Mat_{m\times s}(\Lambda)$ be a coefficient matrix for a system of
defining equations for $W$.  By changing the order of the variables
if necessary, we can use Gaussian elimination on $A(\Psi)$ to obtain
\[
  G(\Psi) = [I_m, C(\Psi)],
\]
where $C(\Psi) \in \Mat_{m \times (s-m)}(\Lambda)$.  Both $A(\Psi)$
and $G(\Psi)$ provide coefficient matrices for equations for $W$,
and so it suffices to show that $C(\Psi)$ has entries in $\ok(t)$.

Let $E/\FF_q(t)$ be a finite extension of fields.  Since $E
\otimes_{\FF_q(t)} W$ is invariant under $\Gamma(E)$, it follows
that, for every $\gamma \in \Gamma(E)$, the matrix
$G(\Psi\gamma^{-1})$ is also the coefficient matrix of a defining
set of equations for $E \otimes_{\FF_q(t)} W$.  Now the columns of
the matrix $[-C(\Psi),\ I_{s-m} ]^{\tr} \in \Mat_{m \times
s}(\Lambda)$ form a basis for $W$.  Thus,
\[
  \bigl[ I_m \ \ C(\Psi\gamma^{-1}) \bigr] \cdot
  \left[ \begin{matrix} -C(\Psi) \\ I_{s-m}
  \end{matrix} \right] = 0, \quad \forall\,\gamma \in
  \Gamma(E),
\]
and so $C(\Psi\gamma) = C(\Psi), \quad \forall\, \gamma \in
\Gamma(E)$.  After varying over all $E/\FF_q(t)$ finite, it follows
from Theorem~\ref{T:GammaFixedIsoK} that $C(\Psi) \in \Mat_{m \times
(s-m)}(\ok(t))$.
\end{proof}

\begin{proposition} \label{P:SubRepsofImages}
Let $M$ be a $t$-motive.  Suppose that $\Phi \in \GL_r(\ok(t))$
represents multiplication by $\bsigma$ on $M$ and that $\Psi \in
\GL_r(\LL)$ is a rigid analytic trivialization for $\Phi$.  For
every $t$-motive $N$ in $\cT_M$ and every $\Gamma$-subrepresentation
$V$ of $N^{\rB}$, there is a sub-$t$-motive $P \subseteq N$ so that
$\xi_M(P) = V$.
\end{proposition}

\begin{proof}
Pick a $\ok(t)$-basis $\bn \in \Mat_{s \times 1}$ for $N$ with
$\bsigma \bn = \Phi_N \bn$, and let $\Psi_N \in \GL_s(\LL)$ be a
rigid analytic trivialization for $\Phi_N$. Let $\bv \in \Mat_{v
\times 1}(N^{\rB})$ be an $\FF_q(t)$-basis for $V$, and extend $\bv$
to a basis $\bu$ of $N^{\rB}$, $\bu = [\bv,\bw]^{\tr}$. By
Lemma~\ref{L:RATBijectiveForSigma}, there is a $H(\Psi) \in
\GL_s(\Sigma)$ so that $\bu = H(\Psi)\cdot \bn$.  We note that
$H(\Psi) = \delta^{-1} \Psi_N^{-1}$ for some $\delta \in
\GL_s(\FF_q(t))$ by Proposition~\ref{P:RATrivialization}(b).

Let $E/\FF_q(t)$ be a finite extension of fields, and let $\gamma
\in \Gamma(E)$. The action of $\gamma$ on $E \otimes_{\FF_q(t)}
N^{\rB}$ is given by the restriction of $\Xi^{(E)}$ as in the proof
of Proposition~\ref{P:xiMFullyFaithful} to $E \otimes_{\FF_q(t)}
N^{\rB}$. Thus,
\[
  \xi^{(E)}(\gamma)(\bu) = H(\Psi\gamma) \bn
  = H(\Psi\gamma)H(\Psi)^{-1}\bu.
\]
Since $V$ is invariant under $\Gamma$, it follows that the upper
right $v \times (s-v)$ block of $H(\Psi\gamma)H(\Psi)^{-1}$ is $0$
for every $\gamma \in \Gamma(E)$.  Let $D(\Psi) \in \Mat_{s \times
(s-v)}(\Lambda)$ be the $s-v$ right-most columns of $H(\Psi)^{-1}$,
and consider the subspace $W \subseteq \Mat_{1 \times s}(\Lambda)$,
\[
  W = \{ \bx \in \Mat_{1\times s}(\Lambda) \mid \bx \cdot
  D(\Psi) = 0 \}.
\]
By our considerations on $H(\Psi)$ at the end of the preceding
paragraph, we see from Lemma~\ref{L:InvariantSubspace} that $W$ has
a set of defining equations over $\ok(t)$.  Thus there is a $C \in
\Mat_{v \times s}(\ok(t))$ of maximal rank so that $C \cdot D(\Psi)
= 0$. Extend $C$ to a matrix $B \in \GL_s(\ok(t))$ such that $C$
forms the top rows of $B$.  Now let $\bn' = B\cdot\bn =
[\bp,\bq]^{\tr}$, with $\bsigma\bn' = \Phi'\bn'$, and let $P$ be the
$\ok(t)$-span of $\bp = C\cdot\bn$. Then
\[
  \bsigma \left[ \begin{matrix} \bp \\ \bq \end{matrix} \right]
  = \bsigma (B\cdot\bn) = \bsigma(B H(\Psi)^{-1}H(\Psi)\bn)
  = \bigl(B\cdot H(\Psi)^{-1} \bigr)^{(-1)} H(\Psi)\cdot B^{-1}\cdot
  \left[ \begin{matrix} \bp \\ \bq \end{matrix} \right].
\]
By construction, the upper right-hand $v \times (s-v)$ block of $B
\cdot H(\Psi)^{-1}$ is $0$.  Thus,
\[
  \bsigma \left[ \begin{matrix} \bp \\ \bq \end{matrix} \right]
  = \left[ \begin{matrix} \Phi_P & 0 \\ * & * \end{matrix} \right]
  \cdot \left[ \begin{matrix} \bp \\ \bq \end{matrix} \right] =
  \Phi'\cdot \left[ \begin{matrix} \bp \\ \bq \end{matrix}
  \right].
\]
Since $\Phi' \in \GL_s(\ok(t))$, it follows that $\Phi_P \in
\GL_v(\ok(t))$.  Thus $P$ is a sub-$t$-motive of $N$. Furthermore,
as $H(\Psi)^{-1} = \Psi_N \delta$, $\delta \in \GL_s(\FF_q(t))$, it
follows that $B\cdot H(\Psi)^{-1}$ is a rigid analytic
trivialization of $\Phi'$.  If we set take $\Psi_P$ to be the upper
left-hand block of $B\cdot H(\Psi)^{-1}$, then $\Psi_P$ is a rigid
analytic trivialization for $\Phi_P$. Moreover, it follows that
$P^{\rB} = V$ by Proposition~\ref{P:RATrivialization}(b).
\end{proof}

\begin{proposition} \label{P:EssentialSurjectivity}
Let $M$ be a $t$-motive.  Suppose that $\Phi \in \GL_r(\ok(t))$
represents multiplication by $\bsigma$ on $M$ and that $\Psi \in
\GL_r(\LL)$ is a rigid analytic trivialization for $\Phi$.  To every
representation $W$ in $\Rep{\Gamma}{\FF_q(t)}$ there is a $t$-motive
$N$ in $\cT_M$ so that $W$ is isomorphic to a subquotient of
$\xi_M(N)$.
\end{proposition}

\begin{proof}
The representation $M^{\rB}$ is faithful by
Corollary~\ref{C:MBFaithful}.  Thus any object in
$\Rep{\Gamma}{\FF_q(t)}$ is isomorphic to a subquotient of a direct
sum of representations of the form $(M^{\rB})^u_v :=
(M^{\rB})^{\otimes u} \otimes ((M^{\rB})^{\vee})^{\otimes v}$. Since
$\xi_M(M^u_v) = (M^u_v)^{\rB} \cong (M^{\rB})^u_v$, the proposition
follows.
\end{proof}

\begin{theorem} \label{T:GammaPsiIsGammaM}
Let $M$ be a $t$-motive.  Suppose that $\Phi \in \GL_r(\ok(t))$
represents multiplication by $\bsigma$ on $M$ and that $\Psi \in
\GL_r(\LL)$ is a rigid analytic trivialization for $\Phi$.  Then the
functor
\[
  \xi_M : \cT_M \to \Rep{\Gamma}{\FF_q(t)}
\]
is an equivalence of Tannakian categories.  Equivalently, the
homomorphism $\pi_M : \Gamma \to \Gamma_M$ is an isomorphism over
$\FF_q(t)$.
\end{theorem}

\begin{proof}
By Propositions~\ref{P:xiMFullyFaithful}
and~\ref{P:SubRepsofImages}, the map $\pi_M$ is faithfully flat
\cite[Prop.\ II.2.21(a)]{delmil82}.  By
Proposition~\ref{P:EssentialSurjectivity}, $\pi_M$ is a closed
immersion \cite[Prop.\ II.2.21(b)]{delmil82}.  Thus $\pi_M$ is an
isomorphism of affine group schemes over $\FF_q(t)$.
\end{proof}

\subsubsection{Remark}
Although we have focused on objects in the category $\cT$, the above
theorem is true (with the same proof) if $M$ is replaced by simply a
rigid analytically trivial pre-$t$-motive.

\section{Galois groups and transcendence}

In this section we first recall the linear independence criterion
introduced in \cite{abp04} by Anderson, Brownawell, and the author,
and one of its applications to $t$-motives.  We then link this
together with our study of the Galois groups of certain $t$-motives,
whose matrices representing multiplication by $\bsigma$ have entries
in $\ok[t]$ and whose fundamental matrices have entries in $\EE$.
These $t$-motives include as a subset rigid analytically trivial
Anderson $t$-motives. In what follows our primary goal will be to
consider the fundamental matrix $\Psi$ associated to such a
$t$-motive $M$  and to equate the transcendence degree over $\ok$ of
$\Psi(\theta)$ and the dimension of the Galois group of~$M$.

\subsection{Linear independence criterion}

\begin{theorem}[{\cite[Thm.~3.1.1]{abp04}}] \label{T:ABPLinearIndependence}
Let $\Phi \in \Mat_r(\ok[t])$ be given such that $\det \Phi =
c(t-\theta)^s$, $c\in \ok^\times$, and suppose that $\psi \in
\Mat_{r \times 1}(\EE)$ satisfies
\[
  \psi^{(-1)} = \Phi\psi.
\]
For every $\rho \in \Mat_{1 \times r}(\ok)$ such that
$\rho\psi(\theta) = 0$, there is a $P \in \Mat_{1 \times r}(\ok[t])$
so that $P(\theta) = \rho$ and $P\psi = 0$.
\end{theorem}

\subsubsection{Connection with solutions of $\sigma$-semilinear
equations} At first glance at the above theorem, the solutions
$\psi$ of the $\sigma$-semilinear equation associated to $\Phi$ are
quite special in that their entries are assumed to be in $\EE$.
However, the following proposition demonstrates that this situation
is not unusual.

\begin{proposition}[{\cite[Prop.~3.1.3]{abp04}}]
\label{P:TateIsEntire} Suppose we are given $\Phi \in
\Mat_r(\ok[t])$ and $\psi \in \Mat_{r \times 1}(\TT)$ so that
\[
  \det \Phi(0) \neq 0, \quad \psi^{(-1)} = \Phi\psi.
\]
Then we necessarily have $\psi \in \Mat_{r\times 1}(\EE)$.
\end{proposition}

\subsubsection{Connection with left $\ktsp$-modules}
The following is a variation on \cite[Prop.~4.4.3]{abp04} with
slightly milder hypotheses.  We do not assume that the representing
matrix $\Phi$ is one directly associated to an Anderson $t$-motive.
 However, we do obtain the same equality of dimensions (with the same
proof).

\begin{proposition}[{\cite[Prop.~4.4.3]{abp04}}]
\label{P:ABPSemilinearEqtMotive} Let $\Phi \in \Mat_r(\ok[t])$ and
$\psi \in \Mat_{r \times 1}(\EE)$ be given as in
Theorem~\ref{T:ABPLinearIndependence}.  Let $N$ be the $\ok[t]$-span
in $\EE$ of the entries of $\psi$, and let $V$ be the $\ok$-span in
$\overline{k_\infty}$ of the entries of $\psi(\theta)$. Then
$\rank_{\ok[t]} N = \dim_{\ok} V$.
\end{proposition}

\begin{proof}
Let $N_1 := \{ P \in \Mat_{1 \times r}(\ok[t]) \mid P\psi = 0 \}$.
We then obtain an exact sequence of $\ok[t]$-modules,
\[
  0 \to N_1 \to \Mat_{1\times r}(\ok[t]) \to N \to 0,
\]
where the second map is given by $P \mapsto P\psi$.  It is easy to
check that this is an exact sequence of left $\ktsp$-modules.  Every
$\ok[t]$-basis for $N_1$ can be extended to a basis of
$\Mat_{1\times r}(\ok[t])$, and so the number of $\ok$-linearly
independent relations of $\ok$-linear dependence among the entries
of $\psi(\theta)$ is at least as great as $\rank_{\ok[t]}N_1$.  Thus
$\rank_{\ok[t]} N \geq \dim_{\ok} V$.  Moreover,
Theorem~\ref{T:ABPLinearIndependence} implies that every
$\ok$-linear relation among the entries of $\psi(\theta)$ lifts to a
$\ok[t]$-linear relation among the entries of $\psi$.  Thus
$\rank_{\ok[t]} N \leq \dim_{\ok} V$.
\end{proof}

\subsection{Dimensions and transcendence degrees}

\subsubsection{Rigid analytic trivializations over $\EE$}
Let $M$ be a $t$-motive.  Suppose that $\Phi \in \GL_r(\ok(t)) \cap
\Mat_r(\ok[t])$ represents multiplication by $\bsigma$ on $M$ and
that $\det \Phi = c(t-\theta)^s$, $c \in \ok^\times$.  An important
observation is that by Propositions~\ref{P:RATrivialization}(c)
and~\ref{P:TateIsEntire}, there is a rigid analytic trivialization
$\Psi$ for $\Phi$ such that $\Psi \in \GL_r(\TT)\cap \Mat_r(\EE)$.

\begin{theorem} \label{T:TrDegDimGalGrp}
Let $M$ be a $t$-motive, and let $\Gamma_M$ be its Galois group.
Suppose that $\Phi \in \GL_r(\ok(t)) \cap \Mat_r(\ok[t])$ represents
multiplication by $\bsigma$ on $M$ and that $\det \Phi =
c(t-\theta)^s$, $c \in \ok^\times$. Let $\Psi$ be a rigid analytic
trivialization of $\Phi$ in $\GL_r(\TT) \cap \Mat_r(\EE)$.  Finally,
let $L$ be the subfield of $\overline{k_\infty}$ generated over
$\ok$ by the entries of $\Psi(\theta)$.  Then
\[
  \trdeg{\ok} L = \dim \Gamma_M.
\]
\end{theorem}

\begin{proof}
By Theorem~\ref{T:GammaPsiIsGammaM}, the groups $\Gamma_M$ and
$\Gamma_\Psi$ are isomorphic.  Moreover, by
Theorem~\ref{T:Smoothness}, their dimension is the same as
$\trdeg{\ok(t)} \Lambda$, where $\Lambda = \ok(t)(\Psi) \subseteq
\LL$. Now let $Q = \ok[\Psi(\theta)] \subseteq L$, and let $S =
\ok(t)[\Psi] \subseteq \Lambda$. Then as rings,
\[
  Q \cong \ok[X_{ij}]/\fa, \quad S \cong \ok(t)[X_{ij}]/\fb,
\]
for ideals $\fa$ and $\fb$.   For $d \geq 1$, let $\ok[X_{ij}]_d$
and $\fa_d$ denote the elements of $\ok[X_{ij}]$ and $\fa$ of total
degree $\leq d$, and let $Q_d \subseteq Q$ correspond to their
quotient. Similarly define $\ok(t)[X_{ij}]_d$, $\fb_d$, and $S_d$.

Fix $d \geq 1$.  Now for any $n \geq 1$, the entries of
$\Psi^{\otimes n}$ comprise all monomials of total degree $n$ in the
$\Psi_{ij}$. If $\psi$ is a column of $\Psi^{\otimes n}$, then
$\psi^{(-1)} = \Phi^{\otimes n} \psi$.  Thus let $\opsi \in \Mat_{N
\times 1}(\EE)$ be the column vector whose entries are the
concatenation of $1$ and each of the columns of $\Psi^{\otimes n}$
for $n \leq d$.  (Here $N = (r^{2d+2}-1)/(r^2-1)$.) Then if $\oPhi
\in \Mat_N(\ok[t]) \cap \GL_N(\ok(t))$ is the block diagonal matrix
\[
  \oPhi := [1] \oplus \Phi^{\oplus r} \oplus \bigl( \Phi^{\otimes 2}
  \bigr)^{\oplus r^2} \oplus \cdots \bigl(\Phi^{\otimes d}
  \bigr)^{\oplus r^d},
\]
it follows that $\opsi^{(-1)} = \oPhi\,\opsi$.  Now it is easy to
see that $Q_d$ is the $\ok$-span of the columns of $\opsi(\theta)$
and that $S_d$ is the $\ok(t)$-span of the columns of $\opsi$.
Since $\oPhi$ and $\opsi$ satisfy the hypotheses for
Proposition~\ref{P:ABPSemilinearEqtMotive}, we see that for all $d
\geq 1$,
\[
  \dim_{\ok} Q_d = \dim_{\ok(t)} S_d.
\]
Thus the homogenizations of $Q$ and $S$ have the same Hilbert series
(see \cite[Ch.~VII, \S 12]{zarsam2}), and so $\trdeg{\ok} L =
\trdeg{\ok(t)} \Lambda$.
\end{proof}

\section{Application to Carlitz logarithms}

\subsection{Carlitz logarithms and $t$-motives}

\subsubsection{The power series $L_\alpha$}
For $\alpha \in \ok^{\times}$ with $\absI{\alpha} <
\absI{\theta}^{q/(q-1)}$, define the power series
\[
L_{\alpha}(t) := \alpha + \sum_{i=1}^\infty \frac{\alpha^{q^i}}{(t -
  \theta^q)(t - \theta^{q^2}) \cdots (t - \theta^{q^i})}.
\]
It is easy to show that $L_\alpha \in \TT$ and that moreover,
$L_\alpha(z)$ converges for all $z \in \KK$ with $\absI{z} <
\absI{\theta}^q$.  By \S \ref{SSS:logC}, we see that
\[
  L_\alpha(\theta) = \log_C(\alpha).
\]
Furthermore, as a power series in $\TT$, $L_\alpha$ also satisfies
the functional equation
\begin{equation} \label{E:LalphaTransformation}
  L_\alpha^{(-1)} = \alpha^{(-1)} + \frac{L_\alpha}{t-\theta}.
\end{equation}

\subsubsection{$t$-motives for Carlitz logarithms}
Fix $\alpha_1, \ldots, \alpha_r \in \ok^{\times}$ with
$\absI{\alpha_i} < \absI{\theta}^{q/(q-1)}$ for $i = 1, \dots, r$.
Set
\[
  \Phi := \Phi(\alpha_1, \dots, \alpha_r)
:= \left[ \begin{matrix}
 t-\theta & 0 & \cdots & 0\\
  \alpha_1^{(-1)}(t-\theta) & 1& \cdots &0\\
  \vdots & \vdots & \ddots & \vdots \\
  \alpha_r^{(-1)}(t-\theta) & 0 & \cdots & 1
\end{matrix} \right] \in \Mat_{r+1}(\ok[t]).
\]
Note that $\Phi$ defines a pre-$t$-motive $X := X(\alpha_1, \dots,
\alpha_r)$ that is an extension of $\one^r$ by the Carlitz motive
$C$:
\[
  0 \to C \to X \to \one^r \to 0.
\]
In spite of the restrictions on $\alpha_1, \dots, \alpha_r$, we will
be able to use the objects $X(\alpha_1,\dots, \alpha_r)$ to
accommodate \emph{all} Carlitz logarithms using
Lemma~\ref{L:LogReduction}.

\begin{proposition}
Let $\alpha_1, \ldots, \alpha_r \in \ok^{\times}$ with
$\absI{\alpha_i} < \absI{\theta}^{q/(q-1)}$ for $i=1, \dots, r$.The
pre-$t$-motive $X = X(\alpha_1, \dots, \alpha_r)$ is a $t$-motive.
\end{proposition}

\begin{proof}
  We prove first that $X$ is rigid analytically trivial and then
  that $X$ is an object in $\cT$.  Define
\[
  \Psi := \Psi(\alpha_1, \dots, \alpha_r) :=
\left[ \begin{matrix} \Omega & 0 &\cdots & 0 \\
  \Omega L_{\alpha_1} & 1& \cdots & 0 \\
  \vdots & \vdots & \ddots & \vdots \\
  \Omega L_{\alpha_r} & 0& \cdots & 1
\end{matrix} \right] \in \GL_{r+1}(\TT)
\]
It is a simple matter to check that $\Psi$ is a rigid analytic
trivialization for $\Phi$ using \eqref{E:LalphaTransformation}.  We
note by Proposition~\ref{P:TateIsEntire} that the entries of $\Psi$
are in $\EE$.

Consider the pre-$t$-motive $C \otimes X$.  We claim that $C \otimes
X$ is in the essential image of the functor $\sM \mapsto M : \cA\cR^I
\to \cR$ of Theorem~\ref{T:AndersonIsogenyFullyFaithful}.  By the
definition of the category $\cT$ in \S\ref{SSS:tMotiveDefinition}, it
will follow that $X$ is a $t$-motive.

Let $\sM := \ok[t]^{r+1}$ with standard $\ok[t]$-basis $m_0, \dots,
m_r$.  Letting $\sm := [m_1, \dots, m_r]^{\tr}$, we give $\sM$ the
structure of a left $\ktsp$-module by setting
\[
  \bsigma \sm := (t-\theta) \Phi\sm.
\]
Now $\sM$ sits in an exact sequence of left $\ktsp$-modules,
\[
  0 \to \sC^{\otimes 2} \to \sM \to \sC^r \to 0,
\]
where $\sC$ is the Carlitz motive in the category of Anderson
$t$-motives of \S\ref{SSS:AndersonCarlitzMotive}.  Since $\sC$ and
$\sC^{\otimes}$ are finitely generated as left
$\ok[\sigma]$-modules, so is $\sM$, and it follows from
\cite[Prop.~4.3.2]{abp04} that $\sM$ is free and finitely generated
as a left $\ok[\bsigma]$-module. Without much difficulty, one shows
that $\bsigma \sM = \langle (t-\theta)^2 m_0, (t-\theta)m_1, \dots,
(t-\theta)m_r \rangle_{\ok[t]}$.   Thus $(t-\theta)^n \sM \subseteq
\bsigma \sM$ for all $n \geq 2$, and $\sM$ is an Anderson $t$-motive
by \S\ref{SSS:AndersontMotiveDefinition}.
\end{proof}

\subsection{The Galois group $\Gamma_X$}

We continue with the notations of the previous section, including
choices of $\alpha_1, \dots, \alpha_r \in \ok^{\times}$ with
$\absI{\alpha_i} < \absI{\theta}^{q/(q-1)}$ for $i=1, \dots, r$.

\subsubsection{The group $G$}
Let $G$ be the algebraic subgroup of $\GL_{r+1}$ over $\FF_q(t)$
such that for all $\FF_q(t)$-algebras $R$,
\[
G(R) = \left\{ \left[ \begin{matrix} * & 0 \\ * & I_r \end{matrix}
  \right] \in \GL_{r+1}(R) \right\}.
\]

\subsubsection{Preliminary calculations}
We claim that $\Gamma_X \subseteq G$.  As in \S\ref{SS:GroupGamma},
we can construct the coordinate ring as the image of $\mu :
\FF_q(t)[X,1/\det X] \to \LL \otimes_{\ok(t)} \LL$, the
$\FF_q(t)$-algebra homomorphism that sends $X$ to $\tPsi =
\Psi_1^{-1}\Psi_2$. As before let $\fq = \ker \mu$. Direct
calculation verifies that $X_{ij} - \delta_{ij} \in \fq$ for all $i
\geq 1$ and $j \geq 2$, where $\delta_{ij}$ is the usual Kronecker
delta.  Thus, $\Gamma_X \subseteq G$. It will be convenient
henceforth to label the non-trivial coordinates of $G \subseteq
\GL_{r+1}$ as $X_0, \dots, X_r$.

Because the Carlitz motive $C$ is contained in $X$, it is an object
in $\cT_X$, and hence there is a surjection $\pi : \Gamma_X
\twoheadrightarrow \Gm$ over $\FF_q(t)$ by
Theorem~\ref{T:GaloisGroupCarlitz}.  Now under $\nu : \ok(t)[X_0,
X_0^{-1}, X_1, \dots, X_r] \to \LL$, which takes $X$ to $\Psi$, we
have $\nu(X_0) = \Omega$.  Thus the action of any $\gamma \in
\Gamma_X(\oFqt)$ on $\Omega$ agrees with the action of the
$X_0$-coordinate of $\gamma$ on $\Omega$.  That is, the surjection
$\pi$ coincides with the natural projection on the $X_0$-coordinate
of $G$.  Let $V$ be the kernel of $\pi$ so that we have an exact
sequence of algebraic groups over $\FF_q(t)$,
\[
  1 \to V \to \Gamma_X \to \Gm \to 1.
\]
The group $V$ is a subgroup of the group of unipotent matrices of
$G$, which itself is naturally isomorphic to $\Ga^r$.  Thus we can
think of $V \subseteq \Ga^r$ with coordinates $X_1, \dots, X_r$.

\begin{proposition}
With notation as above, the group $V$ is a linear subspace of
$\Ga^r$ over $\FF_q(t)$.
\end{proposition}

\begin{proof}
Since $\Gamma_X$ is a smooth over $\FF_q(t)$ by
Theorem~\ref{T:Smoothness}, one verifies that the map $\pi$ is
surjective on Lie algebras, and hence $V$ is also smooth.  Thus it
is determined by the Zariski closure of $V(\oFqt)$ in $\Ga^r$.
Because $\pi$ is surjective, for any non-zero $\alpha \in \oFqt$, we
can choose $\gamma \in \Gamma_X(\oFqt)$ so that $\pi(\gamma) =
\alpha$. Suppose that $\mu = \left[ \begin{smallmatrix} 1 & 0 \\ v &
I_r \end{smallmatrix} \right] \in V(\oFqt)$.  Then direct
calculation gives $\gamma^{-1}\mu\gamma = \left[ \begin{smallmatrix}
1 & 0 \\ \alpha v & I_r \end{smallmatrix} \right] \in V(\oFqt)$, and
thus $V(\oFqt)$ is a linear subspace of $\Ga^r(\oFqt)$.  Since $V$
is smooth, its defining equations over $\oFqt$ are linear forms in
$X_1, \dots, X_r$.  These forms can be defined over $\FF_q(t)$ since
$V$ is simply a linear subspace.
\end{proof}

\subsubsection{Defining polynomials for $\Gamma_X$}
Because the map $\Gamma_X \to \Gm$ is a smooth morphism over
$\FF_q(t)$, Hilbert's Theorem~90 provides an exact sequence
\[
  1 \to V(\FF_q(t)) \to \Gamma_X(\FF_q(t)) \to \Gm(\FF_q(t)) \to 1
\]
by \cite[\S 18.5]{water}.  Let $b_0 \in \FF_q(t)^\times \setminus
\FF_q^{\times}$, and fix a matrix
\begin{equation} \label{E:b0lift}
\gamma = \left[ \begin{matrix} b_0 & 0 & \dots & 0 \\ b_1 & 1 & \cdots & 0 \\
\vdots & \vdots & \ddots & \vdots \\ b_r & 0& \cdots & 1
\end{matrix} \right] \in \Gamma_X(\FF_q(t))
\end{equation}
One checks that the Zariski closure in $\Gamma_X$ of the cyclic
group generated by $\gamma$ is the line in $G$ connecting $\gamma$
to the identity matrix.  Translating this line by any element of $V$
shows that $\Gamma_X$ contains the linear space spanned by $V$ and
$\gamma$.  Since $\Gamma_X$ is irreducible and of dimension $1$
greater than the dimension of $V$, we see that $\Gamma_X$ is this
linear subspace.  Moreover, this implies the following proposition.

\begin{proposition} \label{P:GammaXDefiningEquations}
Suppose $F_1, \dots, F_s \in \FF_q(t)[X_1, \dots, X_r]$ are linear
forms defining $V$, and suppose that $\gamma \in \Gamma_X(\FF_q(t))$
is defined as in \eqref{E:b0lift}.  Then the linear polynomials in
$\FF_q(t)[X_0, \dots, X_r]$,
\[
  G_i := (b_0 - 1)F_i - F_i(b_1, \dots, b_r)(X_0-1), \quad i = 1,
  \dots, s,
\]
are defining polynomials for $\Gamma_X$.
\end{proposition}

\subsection{Linear relations among Carlitz logarithms}

\subsubsection{Defining polynomials for $Z$}
\label{SSS:DefiningPolysForZ} As usual let $Z := \Spec \Sigma$,
where $\Psi = \Psi(\alpha_1, \dots, \alpha_m)$. From
Proposition~\ref{P:ZGammaKProduct} we see that $Z$ and $\Gamma_X$
are isomorphic over $\okt$. Since $\Gamma_X$ is a linear space, $Z$
is also a linear space and isomorphic to $\Gamma_X$ over $\ok(t)$.
Thus we can pick
\[
  \zeta = \left[ \begin{matrix} f_0 & 0 & \dots & 0 \\ f_1 & 1 & \cdots & 0 \\
\vdots & \vdots & \ddots & \vdots \\ f_r & 0& \cdots & 1
\end{matrix} \right] \in Z(\ok(t)),
\]
and then
\[
  Z(\ok(t)) = \zeta\cdot \Gamma_X(\ok(t)).
\]
It is a simple matter to check that the linear polynomials in
$\ok(t)[X_0,\dots,X_r]$,
\[
  H_i := G_i - X_0 G(f_0,\dots,f_r)/f_0, \quad i=1, \dots, s,
\]
are defining polynomials for $Z$.

The following theorems show how the above constructions can be used
to characterize all $k$-linear relations among $\tpi$,
$\log_C(\alpha_1), \dots, \log_C(\alpha_r)$.

\begin{theorem} \label{T:CarlitzLinearRelations}
Let $\alpha_1, \dots, \alpha_r \in \ok^{\times}$ with
$\absI{\alpha_i} < \absI{\theta}^{q/(q-1)}$ for $i=1, \dots, r$. Let
$X= X(\alpha_1, \dots, \alpha_r)$ be the associated $t$-motive.
\begin{enumerate}
\item[(a)] Let $F = c_1X_1 + \dots + c_rX_r$, $c_1, \dots, c_r \in
\FF_q(t)$, be a defining linear form for $V$ so that $G = (b_0-1)F -
F(b_1, \dots, b_r)(X_0-1)$, $b_0, \dots, b_r
  \in \FF_q(t), b_0 \notin \FF_q$, is a defining polynomial for $\Gamma_X$.  Then
\[
  (b_0(\theta)-1)\sum_{i=1}^r c_i(\theta)\log_C(\alpha_i) -
  \sum_{i=1}^r c_i(\theta)b_i(\theta)\tpi = 0.
\]
\item[(b)] Every $k$-linear relation among $\tpi$,
$\log_C(\alpha_1), \dots, \log_C(\alpha_r)$ is a $k$-linear
combination of the relations from part (a).
\item[(c)] Let $N$ be the $k$-linear span of $\tpi$,
$\log_C(\alpha_1), \dots, \log_C(\alpha_r)$.  Then $\dim \Gamma_X =
\dim_k N$.
\end{enumerate}
\end{theorem}

\begin{proof}
Choose $f \in \ok(t)$ as in \S\ref{SSS:DefiningPolysForZ} so that $H
:= G - f X_0$ is a defining polynomial for $Z$.  Then
\begin{equation} \label{E:HEval}
  H(\Omega, \Omega L_{\alpha_1}, \dots, \Omega L_{\alpha_r}) = G(\Omega, \Omega L_{\alpha_1},
  \dots, \Omega L_{\alpha_r}) - f\Omega = 0.
\end{equation}
We see that
\begin{multline*}
f^{(-1)}\Omega^{(-1)} = \sigma G(\Omega, \Omega L_{\alpha_1}, \dots,
\Omega L_{\alpha_r})  =\Omega G \bigl( t-\theta-1,
\alpha_1^{(-1)}(t-\theta), \dots, \alpha_r^{(-1)}(t-\theta) \bigr)
\\+f\Omega - F(b_1, \dots, b_r)\Omega.
\end{multline*}
The first equality is a consequence of \eqref{E:HEval}, and the
second follows from direct computation. Thus
\[
  (t-\theta)f^{(-1)} - f = G\bigl(t - \theta - 1,
  \alpha_1^{(-1)}(t-\theta), \dots, \alpha_r^{(-1)}(t-\theta)\bigr)
  - F(b_1, \dots, b_r).
\]
The right-hand side is a polynomial in $\ok[t]$, so it follows that
$f$ is regular at $t=\theta$.  Indeed if not, then $f^{(-1)}$ must
have a pole at $t = \theta^{(-1)}$, whence $f$ must also have a pole
at $t=\theta^{(-1)}$.  Continuing in this way we see that if $f$ has
a pole at $t=\theta$, then it must have a pole at each
$t=\theta^{(-i)}$, $i\geq 1$, which is not possible.  By a similar
argument we deduce that $f^{(-1)}$ is also regular at $t=\theta$.
Thus we see that
\[
  f(\theta) = -G(-1,0, \dots, 0)|_{t=\theta} + \sum_{i=1}^r
  c_i(\theta)b_i(\theta)
  = - \sum_{i=1}^r c_i(\theta)b_i(\theta).
\]
Equation \eqref{E:HEval} transforms into
\[
  (b_0-1)\sum_{i=1}^r c_i \Omega L_{\alpha_i} - \sum_{i=1}^r c_ib_i
  (\Omega - 1) - f\Omega = 0.
\]
Dividing through by $\Omega$ and evaluating at $t=\theta$, we obtain
part (a).  Part (b) is a consequence of (a) and (c), since
$\Gamma_X$ is a linear space in $G$ over $\FF_q(t)$.  For part (c),
part (a) implies that $\dim_k N \leq \dim \Gamma_X$, since the
defining polynomials for $\Gamma_X$ generate a set of $k$-linear
relations on $\tpi$, $\log_C(\alpha_1), \dots, \log_C(\alpha_r)$ of
dimension $r+1-\dim \Gamma_X$.  However, $\dim_k N \geq \trdeg{\ok}
\ \ok(\tpi, \log_C(\alpha_1), \dots, \log_C(\alpha_r))$ and the
latter quantity is $\dim \Gamma_X$ by
Theorem~\ref{T:TrDegDimGalGrp}.
\end{proof}

\subsubsection{Example}
Let $\zeta_\theta = \sqrt[q-1]{-\theta}$, let $X$ be the $t$-motive
$X(\zeta_\theta)$ of dimension $2$ over $\ok(t)$, and let $\Psi =
\Psi(\zeta_\theta)$. Since $\zeta_\theta$ satisfies
$\fC_t(\zeta_\theta) = \theta \zeta_\theta + \zeta_\theta^q = 0$, we
see that $\zeta_\theta$ is a $t$-torsion point on the Carlitz
module.  Moreover, $\exp_C(\theta \log_C(\zeta_\theta)) = 0$, and
one calculates that
\[
 \log_C(\zeta_\theta) = \frac{\tpi}{\theta}.
\]
Thus $\Gamma_X$ is $1$-dimensional by
Theorem~\ref{T:CarlitzLinearRelations}(c).  If we consider the
function in $\TT$
\[
  \Upsilon := t L_{\zeta_\theta} - \zeta_\theta(t-\theta),
\]
then $\Upsilon^{(-1)} = \Upsilon/(t-\theta)$. Thus $\Upsilon =
f/\Omega$ for some $f \in \FF_q[t]$ by Lemma~\ref{L:OmegaFE}.
Evaluation at $t=\theta$ shows that $f = -1$ identically. Therefore,
$Z_\Psi$ is defined by
\[
 Z_\Psi : \zeta_\theta(t-\theta) X_0 - tX_1 - 1 = 0.
\]
It follows that the defining equation for $\Gamma_X$ is
\[
  \Gamma_X : tX_1 - X_0 + 1 = 0.
\]
In the notation of Theorem~\ref{T:CarlitzLinearRelations}, we have
\begin{gather*}
  F := X_1, \quad b_0 := t+1, \quad b_1 := 1 \\
  G := tX_1 - X_0 + 1, \quad H := G - fX_0,
  \quad f := \zeta_\theta(t-\theta) - 1.
\end{gather*}

\subsection{Algebraic independence of Carlitz logarithms}

Before proving the main result on Carlitz logarithms, we prove a
reduction lemma.

\begin{lemma} \label{L:LogReduction}
Let $\lambda \in \KK^{\times}$. If $\exp_C(\lambda) \in
\ok^{\times}$, then there is an $\alpha \in \ok^{\times}$ with
$\absI{\alpha} < \absI{\theta}^{q/(q-1)}$, an $f \in \FF_q[\theta]$,
and an $n \geq 1$, so that $\lambda = \theta^n\log_C(\alpha) +
f\tpi$.
\end{lemma}

\begin{proof}
Let $\beta = \exp_C(\lambda)$, and assume that $\absI{\beta} \geq
\absI{\theta}^{q/(q-1)}$.  We solve the equation $\fC_t(x) = \theta
x + x^q = \beta$; that is, we find the $t$-division points of
$\beta$ on the Carlitz module.  The Newton polygon for this
equation, along with our assumptions on $\beta$, imply that any
solution $\alpha \in \ok^{\times}$ of this equation must satisfy
$\absI{\alpha} = \absI{\beta}^{1/q}$. Moreover, if for some $\eta
\in \KK$ we have $\exp_C(\eta) = \alpha$, then
\[
  \exp_C(\theta \eta) = \beta = \exp_C(\lambda).
\]
If $\absI{\beta} < \absI{\theta}^{q^2/(q-1)}$, then $\alpha$ is
sufficiently small and we can pick $\eta = \log_C(\alpha)$.  The
result then follows with $n = 1$.  Otherwise, we continue to take
$t$-division values, and for some $n \geq 1$, we have
$\fC_{t^n}(\alpha) = \beta$ with $\absI{\alpha} <
\absI{\theta}^{q/(q-1)}$, for which $\exp_C(\theta^n \log_C(\alpha))
= \beta$.
\end{proof}

\begin{theorem} \label{T:AlgIndCarlitzLogs}
Let $\lambda_1, \dots, \lambda_r \in \KK$ satisfy $\exp_C(\lambda_i)
\in \ok$ for $i=1, \dots, r$.  If $\lambda_1, \dots, \lambda_r$ are
linearly independent over $k$, then they are algebraically
independent over $\ok$.
\end{theorem}

\begin{proof}
Assume that $\lambda_1, \dots, \lambda_r$ are linearly independent
over $k$.  By Lemma~\ref{L:LogReduction}, for each $\lambda_i$ we
can pick $\alpha_i \in \ok^{\times}$ with $\absI{\alpha_i} <
\absI{\theta}^{q/(q-1)}$ so that the $k$-linear span of $\lambda_1,
\dots, \lambda_r$ is contained in the $k$-linear span of $\tpi$,
$\log_C(\alpha_1), \dots, \log_C(\alpha_r)$.  Let $X = X(\alpha_1,
\dots, \alpha_r)$ be the $t$-motive associated to these logarithms
as in the previous sections, and let $\Gamma_X$ be its Galois group.
Let
\[
  L = \ok(\tpi, \log_C(\alpha_1), \dots, \log_C(\alpha_r)),
\]
and let
\[
  N = \textnormal{$k$-linear span of $\tpi$, $\log_C(\alpha_1),
\dots, \log_C(\alpha_r)$}.
\]
Because $\lambda_1, \dots, \lambda_r$ are linearly independent over
$k$, we see that $r \leq \dim_k N \leq r+1$.
Theorems~\ref{T:TrDegDimGalGrp} and~\ref{T:CarlitzLinearRelations}
imply that
\[
  \trdeg{\ok} L = \dim \Gamma_X = \dim_k N.
\]
If $\tpi$, $\log_C(\alpha_1), \dots, \log_C(\alpha_r)$ are linearly
independent over $k$, then they are algebraically independent over
$\ok$, whence the same follows for $\lambda_1, \dots, \lambda_r$
since $L = \ok(\tpi, \lambda_1, \dots, \lambda_r)$. If there is a
linear dependence among $\tpi$, $\log_C(\alpha_1), \dots,
\log_C(\alpha_r)$ over $k$, then $N$ is equal to the $k$-span of
$\lambda_1, \dots, \lambda_r$ and $L = \ok(\lambda_1, \dots,
\lambda_r)$.  Thus in that case $\lambda_1, \dots, \lambda_r$ are
algebraically independent over $\ok$.
\end{proof}

\end{document}